\documentclass[reqno,twoside,11pt,english]{amsart}
\usepackage{amsmath,amsfonts,amssymb,amsthm,epsfig}
\newcommand{\RNum}[1]{\uppercase\expandafter{\romannumeral #1\relax}}
\usepackage{comment}
\usepackage{color}
\usepackage[final,colorlinks,linkcolor=blue,anchorcolor=red,citecolor=blue]{hyperref}
\usepackage{graphicx}
\usepackage{titletoc,epsf}
\usepackage{amsmath,amsfonts,latexsym,amsthm,amsxtra,amssymb,bbm}
\allowdisplaybreaks
\usepackage{graphicx,float}

\newcommand{\lam}{\lambda}
\newcommand{\La}{\Lambda}

\newcommand{\ba}{\backslash}

\newcommand{\la}{\lambda}
\newcommand{\de}{\delta}

\newcommand{\na}{\nabla}

\newcommand{\om}{\omega}
\newcommand{\oo}{\infty}

\newcommand{\non}{\nonumber}

\newcommand{\al}{\alpha}
\newcommand{\dist}{\text{\rm dist}}
\newcommand{\medint}{-\kern -,375cm\int}         
\newcommand{\medintinrigo}{-\kern -,315cm\int}
\newcommand{\wto}{\rightharpoonup}                

\newcommand{\wq}{\infty}
\newcommand{\cal}{\mathcal}

\newcommand{\3}{{\mathcal C}}

\newcommand{\HH}{{\mathcal H}}

\newcommand{\R}{{\mathbb R}}

\newcommand{\N}{{\mathbb N}}

\newcommand{\ep}{\varepsilon}
\newcommand{\su}{\subset}

\newcommand{\ri}{\rightarrow}
\newcommand{\rii}{\rightharpoonup}
\def\supp{\textup{supp}}

\def\loc{\text{\rm loc}}

\newcommand{\sing}{\text{\rm sing}}

\newcommand{\Span}{\text{\rm span}}

\newcommand{\id}{\mathrm{id}}

\newcommand{\e}{\varepsilon}
\newcommand{\Div}{\mathrm{Div}}

\newcommand{\p}{\partial}
\newcommand{\vol}{\mathrm{vol}}
\newcommand{\cJ}{\mathcal{J}}

\numberwithin{equation}{section}
\newtheorem{theorem}{Theorem}[section]

\newtheorem{definition}[theorem]{Definition}

\newtheorem{lemma}[theorem]{Lemma}

\newtheorem{proposition}[theorem]{Proposition}
\newtheorem{remark}[theorem]{Remark}
\newtheorem*{remark*}{Remark}

\theoremstyle{definition}
\newtheorem{Thm}{Theorem}

\makeatother

\begin{document}
	
	\title[Quantitative stratification and optimal regularity for HACS]{Quantitative stratification and optimal regularity for harmonic almost complex structures} 
	
	\author[C.-Y. Guo, M.-L. Liu and C.-L. Xiang]{Chang-Yu Guo, Ming-Lun Liu and Chang-Lin Xiang}
	
	\address[C.-Y. Guo]{Research Center for Mathematics and Interdisciplinary Sciences, Shandong University, 266237, Qingdao, P. R. China, and Department of Physics and Mathematics, University of Eastern Finland, 80101, Joensuu, Finland}
	\email{changyu.guo@sdu.edu.cn}
	
	\address[Ming-Lun Liu]{Research Center for Mathematics and Interdisciplinary Sciences, Shandong University 266237,  Qingdao, P. R. China and and  Frontiers Science Center for Nonlinear Expectations, Ministry of Education, P. R. China}
	\email{minglunliu2021@163.com}

	\address[Chang-Lin Xiang]{Three Gorges Mathematical Research Center, China Three Gorges University, 443002, Yichang,  P. R. China} \email{changlin.xiang@ctgu.edu.cn}

	\thanks{$^*$Corresponding author: Ming-Lun Liu}

	\begin{abstract}
		In a recent interesting work \cite{He-2019}, W.Y. He established the important partial regularity theory and the almost optimal higher regularity theory for energy minimizing harmonic almost complex structures. Based on a new observation on the structure of equations, we give an easier new proof of the partial regularity theorem, and adapting the powerful quantitative stratification method of Naber-Valtorta \cite{Naber-V-2017}, we further prove the rectifiability of singular stratum of energy minimizing harmonic almost complex structures. Based on this, we establish an optimal regularity theory, which improves the corresponding result of He. 
	\end{abstract}
	
	\date{}

	\maketitle
	\tableofcontents
	
	\section{Introduction}
	
	\subsection{Background and motivation}
	In modern differential geometry, almost complex manifolds, which include complex manifolds, symplectic manifolds and K\"ahler manifolds, are one of the central objects with extensive studies.  An almost complex structure, which supports a compatible Riemannian metric, is called an almost Hermitan structure. 
	
	A fundamental question in almost Hermitian geometry, asked by Calabi and Gluck \cite{Calabi-Gluck-1993}, is to find the ``best" almost complex structure among all complex structures that are compatible {with the given almost Hermitian manifold $(M,g)$ of dimension $m=2n$ $(n\ge2)$}. This question was studied in depth by C. Wood in the pioneer works \cite{Wood-1993-Crelle,Wood-1995-Comp} using the theory of twistor bundles. Moreover, motivated by the theory of harmonic mappings, C. Wood introduced in \cite{Wood-1995-Comp} the notion of harmonic almost complex/Hermitian structures, by minimizing the energy functional $E$
	\begin{equation}\label{eq:energy functional}
		E(J)=\int_M|\nabla J|^2d\vol_g,
	\end{equation}
	among all compatible almost complex structures $J$ on $M$. The Euler-Lagrange equation for \eqref{eq:energy functional} is given by 
	\begin{equation}\label{eq:EL by wood}
		[J,\Delta J]=0,
	\end{equation}
	where $\Delta=\Delta_g$ is the rough Laplacian defined by the metric $g$. In certain special cases, harmonic almost complex structures are harmonic mappings into its twistor space and thus in these cases the theory of harmonic maps can be applied directly. However, because of the effect from the underlying metric $g$, one cannot simply apply in parallel the results from harmonic mapping theory in general cases. For instance, harmonic maps with zero energy are simply constant maps with not much interest. However, to determine when there exists a harmonic almost complex structure with zero energy on a given almost Hermitian manifold is a very delicate problem; see \cite{Davidov-2017,He-2019}.

	Denote by $C^\infty(\cJ_g(\Omega)$ the space  of smooth almost complex structures on a domain $\Omega\subset M$:
	$$
	C^\infty(\cJ_g(\Omega)):=\{J\in C^{\infty}(\Omega,T^*M\otimes TM), J^2=-\id,  g(J\cdot,J\cdot)=g(\cdot,\cdot)\}
	$$
	While minimizers of $E$ in $C^\infty(\cJ_g(\Omega)$ receives great attention in the literature (see \cite{Davidov-2017}), from the analytic point view, critical/minimizing points are not always smooth, even if $M$ is Euclidean. Indeed, consider the case $M=\R^4$ and $g=g_0$ is the standard Euclidean metric. Then an energy minimizing harmonic almost complex structure $J$ on $\R^4$ can be viewed as a minimizing harmonic map from $\R^4$ to the sphere $S^2=SO(4)/U(2)$. In this case, $J$ might have singular points. It is thus natural to consider critical/minimizing points of $E$ in $W^{1,2}(\cJ_g(\Omega))$, which consists of admissible almost complex structures on  $\Omega\subset M$:
	$$
	W^{1,2}(\cJ_g(\Omega)):=\{J\in W^{1,2}(\Omega,T^*M\otimes TM), J^2=-\id,  g(J\cdot,J\cdot)=g(\cdot,\cdot)\}.
	$$
	In the context of this situation,  since regularity theory of harmonic mappings is not always applicable, developing a regularity theory for weakly or minimizing harmonic almost complex structures turns out to be indispensable.
	
	The recent interesting work \cite{He-2019} seems to be the first  to explore the regularity theory of minimizing harmonic almost complex structures, where partial regularity theory and almost optimal global regularity were established.
	To record the main results of He \cite{He-2019}, we first introduce the definition of harmonic almost complex structures. 
	\begin{definition}[Harmonic almost complex structures]\label{def:weakly hac}
		An almost complex structure $J\in W^{1, 2}(M^{2n}, T^*M\otimes TM)$ is called a weakly harmonic almost complex structure if it solves \eqref{eq: Harmonic ACS eq} in the weak sense, that is,  for any $T\in  W^{1,2}\cap L^\infty(M, \text{End}(M))$, it holds
		\begin{equation}\label{eq: J weak soultion}
			\int_M \langle \nabla J, \nabla T\rangle d\vol_g+\int_M \langle J\nabla J\cdot \nabla J , T\rangle d\vol_g=0.
		\end{equation}
		
		{An almost complex structure  $J\in W^{1, 2}(M^{2n}, T^*M\otimes TM)$ is called a (locally) minimizing harmonic almost complex structure if it minimizes the energy functional $E$ in $W^{1,2}(\cJ_g(M))$ in the sense that
			$ E(J)\le E(J^\prime) $ for all $J^\prime\in W^{1, 2}(M^{2n}, T^*M\otimes TM)$ such that the support ${\rm supp}(J^\prime-J)$ is compactly contained in $M$. }
	\end{definition}
	
	A starting point of He is to rewrite the Euler-Lagrange equation \eqref{eq:EL by wood} in the form of a semilinear elliptic system:
	\begin{equation}\label{eq: Harmonic ACS eq}
		\Delta J=J\nabla J\cdot\nabla J,
	\end{equation}
	where $J\nabla J\cdot \nabla J$ reads in local coordinates as
	$$
	\left(J\nabla J\cdot \nabla J \right)^i_j=g^{pq}J_j^a(\partial_pJ_a^b)(\partial_qJ_b^i).
	$$
	Based on the new  form \eqref{eq: Harmonic ACS eq}, ideas developed in the theory of harmonic mappings can be borrowed; see \cite{Helein-2002, Lin-Wang-2008-book, Simon-1996-book} for comprehensive introduction on harmonic maps. Recall that a fundamental regularity theory for minimizing harmonic mappings due to R. Schoen and K. Uhlenbeck \cite{SU-1982-JDG}, is the so called $\e$-regularity theorem, which states that if the normalized energy on a ball $B_r(p)\subset M$ is small, then the map is smooth in $B_{r/2}(p)$.  In the context of minimizing harmonic almost complex structures, this basic $\e$-regularity theorem was established by He \cite[Theorem 1]{He-2019}.
	\begin{Thm}[{\cite[Theorem 1]{He-2019}}]\label{thm: ep-regularity}
		Let {$J\in W^{1,2}(\mathcal{J}_g(M))$} be a minimizing harmonic almost complex structure. There exists $\e=\e(m, M, g)>0$ such that for any $p\in M$ if, for some $r\in(0,1)$ we have
		$$
		r^{2-m}\int_{B_r(p)} |\nabla J|^2 d\vol_g<\e,
		$$
		then $J\in C^{\infty}(B_{\frac{r}{2}}(p))$. 
	\end{Thm}
	
	The proof of Theorem \ref{thm: ep-regularity} (see \cite[Section 4]{He-2019}) follows a similar strategy as that of Schoen-Uhlenbeck \cite{SU-1982-JDG}. In particular, one needs to construct comparison almost complex structures and then do some delicate estimates for the comparison structures. This makes the proof somehow quite complicated. 
	
	{ In the present work, we made a further observation on the equation \eqref{eq: Harmonic ACS eq}, that is, we find the following  divergence-free structure (see Section \ref{sec:partial regularity} \eqref{eq:key observation} for more details): }
	$$\Div(J\nabla J)=\nabla J\cdot\nabla J+J\Delta J=J\left(\Delta J-J\nabla J\cdot\nabla J\right)=0.$$
	This allows us to prove a  more general partial regularity result for \emph{weakly harmonic almost complex structures}, Theorem \ref{thm: new ep-regularity}, with a much shorter/easier proof, using the well-known Hardy-BMO duality, which can be viewed as a parallel result  for sphere-valued harmonic mappings established by Evans \cite{Evans-1991}.
	
	 Once a partial regularity as above has been established, as in the theory of harmonic mappings, one can investigate further structure of the singular set
	 $$\sing(J):=\{a\in B_1: J\text{ is not continuous in any neighborhood of }a\},$$ by making use of the tangent harmonic almost complex structures. A tangent almost complex structure can be identified with the energy minimizing harmonic mapping in $W^{1,2}(\R^{2n},SO(2n)/U(n))$ of homogeneous degree zero; see \cite{Davidov-2017,He-2019} for details. For each $k=1,2,\cdots,2n-3$, define $\mathcal{S}^k$ as
	\[
	\mathcal{S}^k=\{x\in \sing(J): \text{no tangent almost complex structure of } J \text{ at } x \text{ is }(k+1)\text{-symmetric} \}.
	\]  
	Then there is a natural stratification of $\sing(J)$:
	\[
	\mathcal{S}^0\subset \mathcal{S}^1\subset \cdots\subset \mathcal{S}^{2n-3}=\sing(J).
	\]
	Moreover, it was shown in \cite[Theorem 1]{He-2019} that for each $1\leq k\leq 2n-3$, $\dim_{\mathcal{H}}\mathcal{S}^k\leq k$. This is parallel to the classical stratification of singular strata of harmonic mappings; see e.g. \cite{Simon-1996-book}.
	
	To study stratification strata of minimizing harmonic mappings, the recent seminal work \cite{Cheeger-Naber-2013-CPAM} of J. Cheeger and A. Naber  introduced a powerful new technique, called \emph{quantitative stratification} in literature, in which the authors  successfully applied it  to obtain almost optimal regularity for minimizing harmonic maps; see also \cite{Cheeger-Naber-2013,Cheeger-N-V-2015,Cheeger-J-N-2021,Edelen-E-2019-TAMS,Fu-Wang-Zhang-2024,GJXZ-2024,Naber-V-2020-JEMS} for other successful applications of this technique.  W.Y. He \cite{He-2019} also successfully adapted this technique to derive almost optimal regularity theory for minimizing harmonic almost complex structures. To record and also to sharpen his result, we need to introduce a few concepts.
	
	{For an almost complex structure  $J\in W^{1,2}(\cJ_g(B_3(p)))$, $r$, $\epsilon>0$ and $k>0$, we may define the \emph{$k$-th quantitative stratum} $S_{\epsilon,r}^k(J)\subset B_3(p)$ (see Definition \ref{def: quantitative stratification} below for precise meaning). If we set $S_{\epsilon}^k(J)=\bigcap\limits_{r>0}S_{\epsilon,r}^k(J)$ and $S^k(J)=\bigcup\limits_{\epsilon>0}S_{\epsilon}^k(J)$, then $S^k(J)$ coincides with the classical $k$-th stratum $\mathcal{S}^k$ defined in the above; see also Lemma \ref{lemma: stratification for SHACS} below.}

	{For each $\La>0$ and $B_1(p)\subset M$, we also need to introduce the space
		$$W^{1,2}_\La(\cJ_g(B_1(p))),$$  which consists of all maps $J\in W^{1,2}(\cJ_g(B_r(p)))$ such that $
		\Phi_{J}(p,1)\le \La,
		$, where $\Phi$ is defined as
		\begin{equation}\label{eq: mono energy}
			\Phi_{J}(p, r)=e^{C_0r}r^{2-m}\int_{B_r(p)}|DJ|^2dx+C_0r.
	\end{equation}}
Here $C_0$ is a constant from the monotonicity formula (see Theorem \ref{thm: monotonicity formula}). 	With these concepts at hand, the second main result of He can be stated as follows. 
	\begin{Thm}[{\cite[Theorem 2]{He-2019}}]\label{thm:B}
		Let {$J\in W^{1,2}_\Lambda(\cJ_g(M))$} be a minimizing harmonic almost complex structure. Then for each $\ep>0$ there exists $C_\ep=C_\ep(m,\Lambda,g,\ep)$ such that for all $r\in(0,1)$ {and $x\in M$,
			\begin{align}\label{eq: Vol estomate}
				\mbox{\rm Vol}(T_r(S^k_{\ep,r}(J))\cap B_1(x))\leq C_\ep r^{m-k-\e}.
		\end{align}}
		Consequently, for any $p\in (2,3)$, there exists a constant $C_0=C_0(m,g,\Lambda,p)>0$ such that
		\[	\max\left\{\int_{B_1(x)}|\nabla J|^p d\vol_g, \int_{B_1(x)}|\nabla^2 J|^{p/2}d\vol_g \right\}\leq \int_{B_1(x)}r_J^{-p}d\vol_g\leq C_0. 
		\]
	\end{Thm}	
	Here the function $r_J$ is called regularity scale function, defined as follows.
	\begin{definition}[Regularity scale]\label{def: regularity scale}
		Let {$J\in W^{1,2}_\Lambda(\cJ_g(M))$} be a minimizing harmonic almost complex structure. For {$x\in B_1(p)$}, we denote the regularity scale $r_J(x)$ of $J$ at $x$ by
		$$r_J(x):=\max\left\{0\leq r\leq 1:\sup_{B_r(x)}\left(r|\nabla J|+r^2|\nabla^2 J|\right)\leq1\right\}.$$
	\end{definition}

	As we have pointed out earlier, in certain special cases, harmonic almost complex structures are simply harmonic mappings into its twistor space. While in the case of minimizing harmonic maps, in the more recent remarkable work, A. Naber and D. Valtorta \cite{Naber-V-2017} were able to prove the rectifiability of each stratum $S^k$ and improve the similar estimate as \eqref{eq: Vol estomate} by removing the extra $\e$ on the exponent of $r$ on the right hand side of \eqref{eq: Vol estomate}. It is thus natural to ask\footnote{This question was initially pointed out by He in the arXiv version Remark 5.13 of \cite{He-2019}. }
	\medskip
	
	\textbf{Question.} \emph{Is it true that each stratum $S^k(J)$ is $k$-rectifiable? Can we remove the extra $\e$ on the exponent of $r$ on the right hand side of \eqref{eq: Vol estomate}?  Namely, whether the following sharper estimate holds for $p\in M$:}
	\begin{align}\label{eq:expected}
		\mbox{\rm Vol}(T_r(S^k_{\ep,r}(u))\cap B_1(p))\leq C_\ep r^{m-k}.
	\end{align}  
	\medskip

	\subsection{Main results}	
The main motivations of this paper is to provide an affirmative answer to the above question, and then based on which we will be able to optimize the regularity result of Theorem \ref{thm:B}.  Consequently, our main results below build a relatively complete and optimal regularity theory for minimizing  harmonic almost complex structure. It is  expected that these results will be useful for further research on almost complex structures. 

Our first main result sharpens the volume estimate of He and verifies the expected rectifiability of singular stratum.
	
	\begin{theorem}[Optimal volume estimates and rectifiability for singular strata]\label{thm: stratification of MHACS}
		Let {$J\in W^{1,2}_\Lambda(\cJ_g(M))$} be a minimizing harmonic almost complex structure. 
		
		{\upshape (i)} For each $\ep>0$ there exists $C_\ep=C_\ep(m,\Lambda,g,\ep)$ such that for all $r\in(0,1)$ {and $p\in M$,
			\begin{align}\label{eq: Vol estomate-1}
				\mbox{\rm Vol}(T_r(S^k_{\ep,r}(J))\cap B_1(p))\leq C_\ep r^{m-k}.
			\end{align}
			Consequently, for all $r\in(0,1)$, we have
			\begin{align}\label{eq: Vol estomate-2}
				\mbox{\rm Vol}(T_r(S^k_{\ep}(J))\cap B_1(p))\leq C_\ep r^{m-k}.
		\end{align}}	
		
			{\upshape (ii)} For each $k$, the set $S^k_\ep(J)$ and $S^k(J)$ are $k$-rectifiable and upper Ahlfors $k$-regular. Moreover, for $\HH^k$-almost every $x\in S^k(J)$, there exists a unique $k$-plane $V^k\su T_xM$ such that every tangent almost complex structure of $J$ at $x$ is $k$-symmetric with respect to $V^k$.
	\end{theorem}
	Note that \eqref{eq: Vol estomate-2} implies that the Minkowski dimension of each $S^k$ satisfies
	\[\dim_{\rm Min} S^k\le k.\]
	In the above theorem, by saying that a subset $A\subset \R^n$ is upper Ahlfors $k$-regular, we mean  there is a constant $M>0$ such that
	\[\cal{H}^k(A\cap B_r(x))\le Mr^k\qquad \text{for all }\,x\in A \text{ and } 0<r<\text{diam}(A). \]
	As a standard application of the improved volume estimate \eqref{eq: Vol estomate-1} in Theorem \ref{thm: stratification of MHACS}, we obtain optimal higher regularity for harmonic almost complex structures, which sharpens the regularity result in Theorem \ref{thm:B}. In below, we denote by $L^q_{weak}(B_1(p))$ the space of weakly $L^q$-integrable functions on $B_1(p)$. 
	
	\begin{theorem}[Regularity estimates on minimizing harmonic almost complex structure]\label{thm: regularity estimates on minimizing HACS}
		There exists a positive constant  $C=C(m,g,\La)$ such that, for any minimizing harmonic almost complex structure $J\in W^{1,2}_\Lambda(\cJ_g(M))$ and $p\in M$, there holds
		\begin{equation}\label{eq: nabla J in L3,infty}
			\begin{aligned}
				{\rm Vol}\Big(\{x\in B_1(p):r|\nabla J|+r^2|\nabla^2 J|>1\}\Big)\leq{\rm Vol}(\{x\in B_1(p):r_J(x)<r\})\leq Cr^3.
			\end{aligned}
		\end{equation}
		In particular,  both $|\na^{\ell} J|$ and $r_J^{-\ell}$ have uniform bounds in {$L^{3/\ell}_{weak}(B_1(p))$} for $\ell=1,2$.
	\end{theorem}
	
	We remark that since $\sing(J)\subset \{x:r_J(x)<r \} $ holds for all $r>0$, the second inequality of \eqref{eq: nabla J in L3,infty}  implies that the singular set of $J$ satisfying  $$\dim_{\rm Min}(\sing(J))\le m-3,$$ which reproduce the Hausdorff dimension estimate of He \cite{He-2019}. We also remark that the estimate of Theorem \ref{thm: regularity estimates on minimizing HACS} is optimal. To see this, consider the case $M=\R^4$ and $g=g_0$ is the standard Euclidean metric. Then an energy minimizing harmonic almost complex structure $J$ on $\R^4$ can be viewed as a minimizing harmonic map from $\R^4$ to the sphere $S^2=SO(4)/U(2)$. In this case, Theorem \ref{thm: regularity estimates on minimizing HACS} reduces to the regularity theorem of Naber-Valtorta \cite{Naber-V-2017}. This shows the sharpness of  Theorem \ref{thm: regularity estimates on minimizing HACS}. 
	
	Now let us explain the strategy we use and the main difficulty we will meet in order to prove the above two theorems. As aforementioned, the main tool is  the robust approach of Naber-Valtorta \cite{Naber-V-2017}.  However, to apply their method, a major technical difficulty we met is the monotonicity formula for minimizing harmonic almost complex structures. More precisely, in Theorem \ref{thm: monotonicity formula} below, we proved the following monotonicity inequality
	\begin{equation}\label{eq: monotonicity formula 1}
		e^{C_0r}r^{2-m}\int_{B_r(a)}|D J|^2dx+\int_{B_R\backslash B_r(a)}\rho^{2-m}|\partial_\rho J|^2dx\le e^{C_0R}R^{2-m}\int_{B_R(a)}|D J|^2dx+C_0(R-r).\tag{MI}
	\end{equation}
	Comparing with that of harmonic maps considered by Naber-Valtorta \cite{Naber-V-2017} or other applications of the Naber-Valtorta theory \cite{Cheeger-N-V-2015,Cheeger-J-N-2021,Edelen-E-2019-TAMS,Fu-Wang-Zhang-2024,GJXZ-2024,Naber-V-2020-JEMS}, in \eqref{eq: monotonicity formula 1}, we do not have an equality, but merely a ``monotonicity inequality with an extra error term". As a result, we have to thoroughly check the whole procedure of \cite{Naber-V-2017}, where a monotonicity equality was used.  It turns out finally that this inequality is sufficient to do the quantitative regularity analysis as that of \cite{Naber-V-2017}. 
	
	Meanwhile, another two important ingredients of the argument that we establish in this paper  are the important compactness property and unique continuation property for minimizing harmonic almost complex structures, which seems to be unavailable in the literature before this paper. We also explore some basic properties of tangent almost complex structures. These results play a fundamental role in adapting the quantitative $\e$-regularity theorem of minimizing harmonic maps to minimizing harmonic almost complex structures. We believe that these results also have independent interest and potential applications in other problems. 
	
	Finally we remark that, in this paper we did not consider the regularity theory for stationary harmonic almost complex structures. The reason is that unlike the case of stationary harmonic maps \cite{Lin-1999-Annals}, the compatibility condition of $J$ with the metric $g$ gives restriction for the class of admissible test functions. Thus it is not clear how to obtain the monotonicity formula from the stationary condition for $J$. It would be interesting to know whether one can actually extend the Naber-Valtorta theory to stationary harmonic almost complex structures. On the other hand, in \cite{He-2023,He-Jiang-2021}, the theory of biharmonic almost complex structures and polyharmonic almost complex structures were studied. Comparing with the case of minimizing biharmonic maps \cite{GJXZ-2024}, it is plausible that one can extend Theorems \ref{thm: stratification of MHACS} and \ref{thm: regularity estimates on minimizing HACS} to biharmonic/polyharmonic almost complex structures. 
	\medskip 
	
	\textbf{Organization of this paper.} This paper is organized as follows. In Section \ref{sec:partial regularity}, we establish the partial regularity theory, compactness theory and unique continuation property for minimizing harmonic almost complex structures. In Section \ref{sec:quantitative ep regularity}, we use the results from Section \ref{sec:partial regularity} to develop the quantitative $\e$-regularity theorem. In Section \ref{sec:reifenberg}, we recall the rectifiable Reifenberg theorems from \cite{Naber-V-2017} and in Section \ref{sec:covering lemma}, we adapt the Naber-Valtorta covering argument to minimizing harmonic almost complex structures. In the final section, Section \ref{sec:proofs}, we prove our main results. 
	
	\medskip
	
	\textbf{Notation of this paper.}  	
	Through the article, we will use the notation $\R^m=\R^{2n}$. The covariant derivative on $M$ is denoted by $\nabla$, while $D, \p$ denote derivatives with respect to local coordinates. Note that $\nabla$ is equivalent to $D$ when the domain of definition is flat. The integration with the metric measure is written as $\int f d\vol_g$, and $\int f:=\int f dx$ for the integration with the Euclidean measure. The symbol $c_n$ denotes a uniformly bounded dimensional constant.

	\section{Partial regularity and compactness theory}\label{sec:partial regularity}
	
	In this section we collect some basic results for minimizing harmonic almost complex structure. Since our discussions are mainly local, after a definite scaling if necessary, we assume that the injectivity radius of $g$ is larger than $10$. Hence for any $p\in M$, the exponential map $\exp_p: B_{10}(0)\subset \R^{m}\rightarrow B_{10}(p)\subset M$ is a diffeomorphism. When we work locally, we can consider the pullback map $\exp^* J=J\circ\exp$ and $\exp^* g=g\circ\exp$ on $B_3=B_3(0)\subset \R^m$, which we will still denote as $g$ and $J$ for brevity.

	\subsection{Monotonicity formula and partial regularity}\label{subsec:mono for and partial regularity}
	
	In this subsection we use Riemannian normal coordinates. We may assume that on $B_1$, $g_{ij}(0)=\delta_{ij}$, $\partial g(0)=0$, and
	\begin{equation}\label{eq: g near Euclidean}
		(1-\delta |x|^2)\delta_{ij}\leq g(x)\leq (1+\delta |x|^2)\delta_{ij}, |\p g|(x)\leq \delta |x|, |\p^2 g|\leq \delta
	\end{equation}
	for some small but fixed constant $\delta$. 
	
	{
		Since the almost complex structure $J$ carries the geometric information of $M$, we cannot always reduce $g$ to the Euclidean metric $\delta_{ij}$. Thus it is useful to compare the energy measured by $g$ with that by Euclidean metric $\delta_{ij}$, for which we state as follows.} 
	
	\begin{proposition}[{\cite[Proposition  3.1]{He-2019}}]\label{prop: comparison}
		The following comparison result holds: 
		\begin{equation}\label{eq: comparison1}
			(1-c_n \delta r^2)\int_{B_r} |D J|^2 -c_n\delta^2 r^m\le\int_{B_r} |\nabla J|^2d\vol_g\leq (1+c_n \delta r^2)\int_{B_r} |D J|^2 +c_n\delta^2 r^m.
		\end{equation}
		Consequencely, for $c_n\delta r^2\le\frac12$, it holds  
		\begin{equation}\label{eq: comparison2}
			(1-c_n \delta r^2)\int_{B_r} |\nabla J|^2d\vol_g -c_n\delta^2 r^m\le\int_{B_r} |D J|^2\leq (1+2c_n \delta r^2)\int_{B_r} |\nabla J|^2d\vol_g +2c_n\delta^2 r^m.
		\end{equation}
	\end{proposition}

	Next, we derive an improved monotonicity formula based on the result of He \cite[Theorem 3.1]{He-2019}.
	{\begin{theorem}[Monotonicity formula]\label{thm: monotonicity formula}
			{Let $J\in W^{1,2}(\mathcal{J}_g(B_1(p)))$ be a minimizing harmonic almost complex structure on $B_1(p)\subset M$.} Then, for any $a\in B_{\frac{1}{2}}(p)$ and every $0 < r< R \le\frac12$, there holds
			\begin{equation}\label{eq: monotonicity formula}
				e^{C_0r}r^{2-m}\int_{B_r(a)}|D J|^2dx+\int_{B_R\backslash B_r(a)}\rho^{2-m}|\partial_\rho J|^2dx\le e^{C_0R}R^{2-m}\int_{B_R(a)}|D J|^2dx+C_0(R-r).
			\end{equation}
			where $C_0>0$ is a constant, $\rho=|x-a|$ and $\partial_{\rho}J=DJ\cdot\frac{x-a}{|x-a|}$.
	\end{theorem}}
	
	\begin{proof}
		Identify $B_{\frac12}(a)\subset M$ with $B_{\frac12}(0)\subset \R^m$ via the exponential map $\exp$, and denote the pullback maps $\exp_a^*J$ and $\exp_a^*g$ by $J$ and $g$, respectively. For any $0<r\leq \frac12$, consider
		$$
		J_r(x)=\begin{cases}
			J(x), |x| \geq r;\\
			J(\frac{rx}{|x|}), |x|<r.
		\end{cases}$$
		Clearly $J_r\in W^{1,2}\cap L^\infty$, but $J_r$ may not be $g(x)$-compatible. 
		By \cite[Proposition 3.3]{He-2019}, we can construct a $\bar{J_r}$ compatible with $g(x)$ out of $J_r$. By \cite[Proposition 3.4]{He-2019}, we have
		\begin{equation}\label{eq: comparison3}
			(1-c_n \delta r^2)\int_{B_r} |\nabla \bar J_r|^2d\vol_g-c_n\delta^2 r^m \le\int_{B_r} |\nabla J_r|^2d\vol_g\leq (1+c_n \delta r^2)\int_{B_r} |\nabla \bar J_r|^2 d\vol_g+c_n\delta^2 r^m.
		\end{equation}
		
		Since $J$ is minimizing, using \eqref{eq: comparison3} we have$$
		\begin{aligned}
			\int_{B_r}|\nabla J|^2d\vol_g
			\le& \int_{B_r}|\nabla \bar{J_r}|^2d\vol_g\le (1+2c_n \delta r^2)\int_{B_r} |\nabla J_r|^2 d\vol_g+2c_n\delta^2 r^m\\
			\le& (1+4c_n \delta r^2)\int_{B_r} |D J_r|^2 +4c_n\delta^2 r^m\\
			=& (1+4c_n \delta r^2)\int_0^rd\rho\int_{\partial B_\rho} |D J_r|^2(x)d\mathcal{H}^{m-1}(x) +4c_n\delta^2 r^m.
		\end{aligned}
		$$
		Here $$
		\begin{aligned}
			\int_{\partial B_\rho} |D J_r|^2(x)d\mathcal{H}^{m-1}(x)
			=& \int_{\partial B_r} |D J_r|^2\left(\frac{\rho y}{r}\right)\left(\frac{\rho}{r}\right)^{m-1}d\mathcal{H}^{m-1}(y)\\
			=& \int_{\partial B_r} \left(|D J(y)|^2-|\partial_rJ(y)|^2\right)\left(\frac{\rho}{r}\right)^{m-3}d\mathcal{H}^{m-1}(y),
		\end{aligned}
		$$
		which implies
		$$
		\begin{aligned}
			&\int_{B_r}|\nabla J|^2d\vol_g\\
			\le& (1+4c_n \delta r^2)r^{3-m}\left(\int_0^r\rho^{m-3}d\rho\right)\int_{\partial B_r} \left(|D J|^2-|\partial_rJ|^2\right)d\mathcal{H}^{m-1} +4c_n\delta^2 r^m\\
			=& (1+4c_n \delta r^2)\frac{r}{m-2}\int_{\partial B_r} \left(|D J|^2-|\partial_rJ|^2\right) d\mathcal{H}^{m-1} +4c_n\delta^2 r^m.
		\end{aligned}
		$$
		Hence$$
		\begin{aligned}
			&\int_{\partial B_r}|\partial_rJ|^2 d\mathcal{H}^{m-1}\\
			\le& \int_{\partial B_r}|D J|^2 d\mathcal{H}^{m-1}+8mc_n\delta^2r^{m-1}-\frac{m-2}{r}(1-4c_n\delta r^2)\int_{B_r}|\nabla J|^2d\vol_g\\
			\le& \int_{\partial B_r}|D J|^2 d\mathcal{H}^{m-1}+8mc_n\delta^2r^{m-1}-\frac{m-2}{r}(1-4c_n\delta r^2)\int_{B_r}|\nabla J|^2d\vol_g\\
			\le& \int_{\partial B_r}|D J|^2 d\mathcal{H}^{m-1}+10mc_n\delta^2r^{m-1}-\frac{m-2}{r}(1-6c_n\delta r^2)\int_{B_r}|D J|^2,
		\end{aligned}
		$$
		where in the last inequality we used \eqref{eq: comparison1}. Multiply $e^{6mc_n\delta r}r^{2-m}$ on both sides we have$$
		\begin{aligned}
			&r^{2-m}\int_{\partial B_r}|\partial_rJ|^2 d\mathcal{H}^{m-1}\le e^{6mc_n\delta r}r^{2-m}\int_{\partial B_r}|\partial_rJ|^2 d\mathcal{H}^{m-1}\\
			\le& e^{6mc_n\delta r}r^{2-m}\int_{\partial B_r}|D J|^2 d\mathcal{H}^{m-1}\\
			&-e^{6mc_n\delta r}r^{2-m}\frac{m-2}{r}(1-6c_n\delta r^2)\int_{B_r}|D J|^2+e^{6mc_n\delta r}10mc_n\delta^2r\\
			\le& e^{6mc_n\delta r}r^{1-m}\left(r\int_{\partial B_r}|D J|^2 d\mathcal{H}^{m-1}+ 6mc_n\delta r \int_{B_r}|D J|^2+ (2-m)\int_{B_r}|D J|^2\right)\\
			& +e^{6mc_n\delta}10mc_n\delta^2\\
			\le& \frac{d}{dr}\left(e^{C_0r}r^{2-m}\int_{B_r}|D J|^2+C_0r\right),
		\end{aligned}
		$$where $C_0=\max\{6mc_n\delta, e^{6mc_n\delta}10mc_n\delta\}$. Integrating it from $r$ to $R$ gives \eqref{eq: monotonicity formula}.
	\end{proof}
	
	\begin{remark}\label{rmk:on monotonicity formular Euclidean}
		When $M$ is Euclidean, $C_0=0$. Thus, \eqref{eq: monotonicity formula} reduces to the simpler form
		\begin{equation}\label{eq: delta-monotonicity formula}
			r^{2-m}\int_{B_r(a)}|DJ|^2dx+\int_{B_R\backslash B_r(a)}\rho^{2-m}|\partial_\rho J|^2dx\le R^{2-m}\int_{B_R(a)}|DJ|^2dx.
		\end{equation}
	\end{remark}

	For $r\in(0,1)$, recall from \eqref{eq: mono energy}
	$$
	\Phi_{J}(p, r)=e^{C_0r}r^{2-m}\int_{B_r(p)}|DJ|^2dx+C_0r.
	$$
	For our later application, we introduce
	\begin{equation}\label{eq: nablaJ energy}
		\Phi^g_{J}(p, r)=r^{2-m}\int_{B_r(p)}|\nabla J|^2d\vol_g,
	\end{equation}
	and
	\begin{equation}\label{eq: DJ energy}
		\Phi^\delta_{J}(p, r)=r^{2-m}\int_{B_r(p)}|DJ|^2dx.
	\end{equation}
	By the monotonicity formula \eqref{eq: monotonicity formula}, we know that $\Phi_J(x,r)$ is increasing in $r$, and hence $\lim_{r\downarrow 0}\Phi_J(x,r)$ exists.
	\begin{definition}[Density function]\label{def:density function}
		For $x\in M$, we define the density function as
		$$
		\Phi_J(x):=\lim\limits_{r\downarrow0}\Phi_J(x,r)=\lim\limits_{r\downarrow0}\Phi^g_J(x,r)=\lim\limits_{r\downarrow0}\Phi^\delta_J(x,r).
		$$
	\end{definition}

	To state our partial regularity theorem below, we first recall the definition of Morrey space. Let $1\le p<\infty$ and $0\le s\le n$. The Morrey space $M^{p,s}(U)$
	consists of functions $u\in L^{p}(U)$ such that
	$$
	\|u\|_{M^{p,s}(U)}\equiv\sup_{x\in U,0<r<\mathrm{diam}(U)}r^{-s/p}\|u\|_{L^{p}(B_{r}(x)\cap U)}<\infty.
	$$
	The space $M_{1}^{p,s}(U)$ consists of functions in $M^{p,s}(U)$ whose weak gradient belongs to $M^{p,s}(U)$. 
	
	Our main partial regularity theorem reads as follows. In case of harmonic maps, similar result has been  obtained by Evans \cite{Evans-1991}. 
	\begin{theorem}\label{thm: new ep-regularity}
		Let $J\in W^{1,2}(\mathcal{J}_g(M))$ be a weakly harmonic almost complex structure. There exists $\e=\e(m, M, g)>0$ such that if $\|\nabla J\|_{M^{m-2,2}(B_1(p))}<\epsilon$, then $J$ is smooth in $B_1(p)$. 
	\end{theorem}
	
	The followin Hardy-BMO inequality, which is an application of the Hardy-BMO duality and div-curl lemma, is well-known; see e.g. \cite{CLMS-1993,Evans-1991,Schikorra-2010} or \cite[Lemma 2.4]{Guo-Liu-Xiang-2025}.

	\begin{lemma}[Hardy-BMO inequality]\label{lemma: Hardy-BMO duality by L^p-L^q-Morrey}
		For any $p\in (1,\infty)$ and $\alpha\in (1,n)$, there exists a constant $C=C(n,p,\alpha)>0$ such that the following holds: 
		
		{\upshape (1)}. For all balls $B_r(x_0)\subset \R^n$, functions $a\in M_1^{\alpha,n-\alpha}(B_{2r}(x_0))$, $\Gamma\in L^q(B_{r}(x_0),\R^n)$, $b\in W_0^{1,p}\cap L^\infty(B_r(x_0))$ with $\frac{1}{p}+\frac{1}{q}=1$ and $\Div(\Gamma)=0$ in the weak sense on $B_r(x_0)$, we have 
		$$
		\left|\int_{B_r(x_0)}\langle\nabla a,\Gamma\rangle bdx\right|\le C\left\|\Gamma\right\|_{L^q(B_r(x_0))}\|\nabla b\|_{L^p(B_r(x_0))}\|\nabla a\|_{M^{\alpha,n-\alpha}(B_{2r}(x_0))}.
		$$
		
		{\upshape (2)}. For all balls $B_r(x_0)\subset \R^n$, functions $\varphi\in C_0^\infty(B_r(x_0))$, $\Gamma\in L^q(B_{r}(x_0),\R^n)$, $b\in M_{1}^{\alpha,n-\alpha}\cap L^\infty(B_{2r}(x_0))$ with $\frac{1}{p}+\frac{1}{q}=1$ and $\Div(\Gamma)=0$ in the weak sense on $B_r(x_0)$, we have
		$$
		\left|\int_{B_r(x_0)}\langle\nabla \varphi,\Gamma\rangle bdx\right|\le C\left\|\Gamma\right\|_{L^q(B_r(x_0))}\|\nabla \varphi\|_{L^p(B_r(x_0))}\|\nabla b\|_{M^{\alpha,n-\alpha}(B_{2r}(x_0))}.
		$$
	\end{lemma}
	
	Now we are ready to prove Theorem \ref{thm: new ep-regularity}.
	\begin{proof}[Proof of Theorem \ref{thm: new ep-regularity}]
		Fix any $x_0\in B_1(p)\subset M$ and write $B_r$ for the ball {$\exp^{-1}(B_r(x_0))\subset \exp^{-1}(B_1(p))\subset \R^m$}, where $r$ is small enough such that $B_{2r}(x_0)\subset B_1(p)$. Assume $\e\in(0,\e_0)$ for some $\e_0=\e_0(m)>0$ to be determined later. 
		
		{Note that by the definition of $W^{1,2}(\cJ_g)$, we have $J\in L^{\infty}$ since $$
			|J|^2=g^{ij}g_{ks}J_i^kJ_j^s=g^{ij}g(J_i^k\partial_k,J_j^s\partial_s)=g^{ij}g_{ij}=2n.	
			$$}
		By \eqref{eq: comparison2} and our smallness assumption, we  have $$
		\|J\nabla J\|_{M^{2,m-2}(B_r)}= C(m)\|\nabla J\|_{M^{2,m-2}(B_r)}\le C(m)\e.
		$$	
		Since $J^2=-\id$, left multiply both sides of \eqref{eq: Harmonic ACS eq} by $J$, and we get 
		\begin{equation}\label{eq:key observation}
			\Div(J\nabla J)=\nabla J\cdot\nabla J+J\Delta J=J\left(\Delta J-J\nabla J\cdot\nabla J\right)=0.
		\end{equation}
		In order to estimate $\|\nabla J\|_{L^{2}(B_r)}$, we define 
		$$T:=\left\{ \varphi\in C_0^\infty(B_r,\R^{m\times m}):\|\nabla \varphi\|_{L^{2}(B_r)}\le1\right\}.$$
		Then
		\begin{equation}\label{eq: nabla J L2 estimate-1}
			\begin{aligned}
				\|\nabla J\|_{L^2(B_r)}\lesssim& \sup_{\varphi\in T}\left|\int_{B_r}\left\langle\nabla J,\nabla\varphi\right\rangle\right|\lesssim \sup_{\varphi\in T}\left|\int_{B_r}\left\langle\Delta J,\varphi\right\rangle\right|\\
				\lesssim& \sup_{\varphi\in T}\left|\int_{B_r}\left\langle J\nabla J\cdot\nabla J,\varphi\right\rangle\right|
			\end{aligned}
		\end{equation}
		Write $J=(J_i^j)$ and $\varphi=(\varphi_i^j)$. Observe that $$
		\left\langle J\nabla J\cdot\nabla J,\varphi\right\rangle= \left\langle J_i^k\nabla J_k^s,\nabla J_s^j \right\rangle\varphi_i^j.
		$$
		Since $\Div(J_i^k\nabla J_k^s)=0$, we may apply Lemma \ref{lemma: Hardy-BMO duality by L^p-L^q-Morrey} (1) with $\Gamma=J_i^k\nabla J_k^s\in L^2$, $a=J_s^j\in M_1^{\lambda,m-\lambda}$ and $b=\varphi_i^j\in W^{1,2}_0\cap L^\infty$ to obtain
		\begin{equation}\label{eq: nabla J L2 estimate-2}
			\begin{aligned}
				\int_{B_r}\left\langle J\nabla J\cdot\nabla J,\varphi\right\rangle\lesssim& \|J\nabla J\|_{L^2(B_r)}\|\nabla J\|_{M^{2,m-2}(B_{2r})}\|\nabla\varphi\|_{L^2(B_r)}\\
				\lesssim& r^{\frac{m-2}{2}}\e\|\nabla J\|_{M^{2,m-2}(B_{2r})}\|\nabla\varphi\|_{L^2(B_r)},
			\end{aligned}
		\end{equation}
		where in the second inequality we used $J\in L^\infty$ and $$
		\|J\nabla J\|_{L^2(B_r)}\lesssim r^{\frac{m-2}{2}}\|\nabla J\|_{M^{2,m-2}(B_{r})}\lesssim \e r^{\frac{m-2}{2}}.
		$$
		Combining \eqref{eq: nabla J L2 estimate-1} and \eqref{eq: nabla J L2 estimate-2}, we obtain that for some $C_0(m)>0$
		\begin{equation}\label{eq: nabla J L2 estimate-3}
			\|\nabla J\|_{L^2(B_r)}\le C_0\e r^{\frac{m-2}{2}}\|\nabla J\|_{M^{2,m-2}(B_{2r})}.
		\end{equation}
		Letting $\e=\frac{1}{2C_0}$ in \eqref{eq: nabla J L2 estimate-3} gives  
		\begin{equation}\label{eq: nabla J L2 estimate-4}
			\|\nabla J\|_{L^2(B_r)}\le \frac12 r^{\frac{m-2}{2}}\|\nabla J\|_{M^{2,m-2}(B_{2r})}.
		\end{equation}
		Clearly we can replace $B_{2r}(x_0)$ with any $B_s(y_0)\subset B^m_1(0)$ containing $B_{2r}(x_0)$.
		Thus \eqref{eq: nabla J L2 estimate-4} implies 
		\begin{equation*}\label{eq: nabla J L2 estimate-5}
			r^{-\frac{m-2}{2}}\|\nabla J\|_{L^2(B_r(x_0))}\le \|\nabla J\|_{M^{2,m-2}(B_{s}(y_0))},
		\end{equation*}
		which is valid for all $r,s,x_0,y_0$ such that $B_{2r}(x_0)\subset B_s(y_0)\subset  B_1(0)$. Note that the family of balls $\{B_{r}(x_0)\}$ forms an open cover of $B_{\frac{s}{2}}(y_0)$. Thus we can take the supremum over all admissible $B_r(x_0)$ to find
		\begin{equation}\label{eq: nabla J Morrey estimate-1}
			\|\nabla J\|_{M^{2,m-2}(B_{\frac{ s}{2}}(y_0))}\le \frac12\|\nabla J\|_{M^{2,m-2}(B_{s}(y_0))}.
		\end{equation}
		Iterating \eqref{eq: nabla J Morrey estimate-1}, we obtain 
		\begin{equation*}\label{eq: nabla J Morrey estimate-2}
			\|J\|_{M^{2,m-2}(B_{2^{-k}s}(y_0))}\le 2^{-k}\|\nabla J\|_{M^{2,m-2}(B_{s}(y_0))}\quad \text{for all } k\in\N. 
		\end{equation*}
		For $r\approx2^{-k}s$, we have $2^{-k}\approx (r/s)$. This implies 
		$$
		\|J\|_{M^{2,m-2}(B_{r}(y_0))}\le Crs^{-1}\|\nabla J\|_{M^{2,m-2}(B_{s}(y_0))}.
		$$
		We may choose $s=r_0$ fixed such that $r_0^{2-m}\int_{B_{r_0}(y_0)}|\nabla J|^2dx\le \e$. Then for all $r\leq r_0/2$, there holds
		$$
		\|J\|_{M^{2,m-2}(B_{r}(y_0))}\le Crr_0^{-1}\|\nabla J\|_{M^{2,m-2}(B_{r_0}(y_0))}.
		$$
		This gives 
		\begin{equation*}\label{eq: nabla J new Morrey space}
			\nabla J\in M_{\loc}^{2,m-2+2}(B_{\frac{r_0}{2}}(y_0)).
		\end{equation*}
		
		Finally, using Morrey's Dirichlet growth theorem (see for instance \cite{Giaquinta-Book}) we infer that  $J\in C_{\loc}^{0,1/2}(B_{\frac{r_0}{2}}(y_0))$. Then using standard Schauder estimates, we obtain the smoothness of $J$. 
	\end{proof}
	\begin{remark}\label{rmk: 2-dim HACS is smooth}
		Consider the following equation:
		\begin{equation}\label{eq: B^2 to SO(2n)/U(n)}
			\Delta J=J\na J\cdot\na J, \qquad J\in W^{1,2}(B_1^2, TM\otimes T^*M).
		\end{equation} 
		Let $J=\left(J_1,J_2,...,J_{2n}\right)$, we can rewrite the equation as
		$$
		\Delta J_i=(J\na J)\cdot\na J_i,\qquad i=1,...,2n.
		$$
		Since $J^2=-\id$ and $J^TJ=\id$, we have $J\nabla J$ is anti-symmetric, and so $J_i$ is continuous by \cite[Theorem I.1]{Riviere-2007-Invent}.
		Then the standard Schauder method implies that the weak solutions of \eqref{eq: B^2 to SO(2n)/U(n)} are smooth.
	\end{remark}

	\subsection{Compactness property}
	
	In this subsection, we prove a compactness theorem for minimizing harmonic almost complex structure, which partially answer \cite[Problem 6.8]{He-2019}.

	\begin{theorem}[Compactness of minimizing harmonic almost complex structure]\label{thm: compact minimizing}
		Let {$J_i\in W^{1,2}(\cJ_g(M))$} be a sequence of minimizing harmonic almost complex structure. If $J_i\wto J$ weakly in $W^{1,2}$, then {$J_i\in W^{1,2}(\cJ_g(M))$} is a minimizing harmonic almost complex structure and $J_i\to J$ strongly in $W^{1,2}$.
	\end{theorem}
	
	In order to establish the compactness property of minimizing  harmonic almost complex structures, we need an extension lemma, which was first partially proved by Schoen-Uhlenbeck \cite{SU-1982-JDG} and later in the following form by Luckhaus \cite{Luckhaus-1988-Indiana}. 
	
	\begin{lemma}\label{lemma: Luckhaus}
		For $n\ge 2$, suppose $u, v\in H^{1}( S^{n-1}, N)$. Then for any $\epsilon\in( 0, 1)$, there is $w\in H^{1}(S^{n- 1}\times [1-\epsilon, 1], \R^L)$ such that $w|_{S^{n-1}\times\{1\}}=u$, $w| _{S^{n-1}\times\{1-\epsilon\}}= v$,
		\begin{equation}\label{eq: Luckhuas Lemma -1}
			\int_{S^{n-1}\times[1-\epsilon,1]}|\nabla w|^{2}\:dx\leq C\epsilon\int_{S^{n-1}}\left(|\nabla_{T}u|^{2}+|\nabla_{T}v|^{2}\right)+C\epsilon^{-1}\int_{S^{n-1}}|u-v|^{2},
		\end{equation}
		and
		\begin{equation}\label{eq: Luckhuas Lemma -2}
			\dist^{2}(w(x),N) \leq C\epsilon^{1-n}\left(\int_{S^{n-1}}\left(|\nabla_{T}u|^{2}+|\nabla_{T}v|^{2}\right)\right)^{\frac{1}{2}} \left(\int_{S^{n-1}}|u-v|^{2}\right)^{\frac{1}{2}}+C\epsilon^{-n}\int_{S^{n-1}}|u-v|^{2}
		\end{equation}
	\end{lemma}
	
	Now we are ready to prove theorem \ref{thm: compact minimizing}.
	
	\begin{proof}[Proof of Theorem \ref{thm: compact minimizing}]
		For any unit ball $B_1\subset\subset M$ and a small $\lambda\in(0,1)$, let $K\in W^{1,2}(M, T^*M\otimes TM)$ be such that $K=J$ on $B_1\backslash B_{1-\la}$. By Rellich's compactness theorem, $J_i\to J$ strongly in $L^2$. Then by Fubini's theorem and Fatou's lemma, there is $\rho\in(1-\la_0, 1)$ such that
		$$
		\lim\limits_{i\to\infty}\int_{\partial B_\rho}|J_i-J|^2d\mathcal{H}^{n-1}=0,\quad
		\int_{\partial B_\rho}\left(|\nabla J_i|^2+|\nabla K|^2\right)d\mathcal{H}^{n-1}\leq C<+\infty.
		$$
		Applying Lemma \ref{lemma: Luckhaus} to $J_i$ and $K$, we conclude that there is $P_i$ such that for suitable $\la_i\downarrow0$, $$
		P_i(x)=
		\begin{cases}
			K\left(\frac{x}{1-\lambda_i}\right), & |x|\leq\rho(1-\lambda_i) \\
			J_i(x), & |x|=\rho 
		\end{cases}
		$$and$$
		\la_i^{-1}\int_{\partial B_\rho}|J_i-K|^2d\mathcal{H}^{n-1}\to 0\text{ as }i\to\infty.$$
		Thus, 
		\begin{equation}\label{eq: estimate for Pi}
			\begin{aligned}
				\int_{B_\rho\setminus B_{\rho(1-\lambda_i)}}|\nabla P_i|^2 &\leq C\lambda_i\int_{\partial B_\rho}\left(|\nabla J_i|^2+|\nabla K|^2\right)d\mathcal{H}^{n-1}+C\lambda_i^{-1}\int_{\partial B_\rho}|J_i-K|^2d\mathcal{H}^{n-1}\\&\to0\text{ as }i\to\infty.
			\end{aligned}
		\end{equation}
		
		Then by minimality and \eqref{eq: estimate for Pi}, we have
		$$
		\begin{aligned}
			\int_{B_\rho}|\nabla J|^2&\leq \quad\lim_{i\to\infty}\int_{B_\rho}|\nabla J_i|^2\leq\quad\lim_{i\to\infty}\int_{B_\rho}|\nabla P_i|^2\\
			&=\quad\lim_{i\to\infty}\left[\int_{B_{\rho(1-\lambda_i)}}\left|\nabla K\left(\frac{x}{1-\lambda_i}\right)\right|^2dx+\int_{B_\rho\setminus B_{\rho(1-\lambda_i)}}|\nabla P_i|^2dx\right]\\
			&\leq\quad\lim_{i\to\infty}\left[(1-\lambda_i)^{n-2}\int_{B_\rho}|\nabla K(y)|^2dy\right]\\
			&\leq\quad\int_{B_\rho}|\nabla K|^2.
		\end{aligned}
		$$
		This implies both minimality of $J$ and strong convergence of $J_i$ to $J$ by letting $K=J$.	
		\end{proof}
		
		\subsection{Unique continuation property}

		The main result of this subsection is the following unique continuation property for minimizing harmonic almost complex structures. 
		\begin{proposition}[Unique continuation property]\label{prop: unique continuation} 
			{Suppose that $J\in W^{1,2}_{\Lambda}(\cJ_g(B_1(p)))$ is a minimizing harmonic almost complex structure on $B_1(p)\subset M$.} Then $J$ enjoys the unique continuation property in the sense that if there is another weakly harmonic almost complex structure $J^\prime$ in $B_{1}(p)$ such that $J=J^\prime$ almost everywhere on an open set, and $J^\prime$ is smooth away from a set $\Sigma'$ of finite $\mathcal{H}^{m-2}$-measure, then $J\equiv J^\prime$ on $B_{1}(p)$. 
		\end{proposition}

		For the proof of Proposition \ref{prop: unique continuation}, we need the following unique continuation theorem for smooth harmonic almost complex structures.
		\begin{theorem}\label{thm: unique continuation}
			Let {$J,J'\in C^\infty(\cJ_g(M))$} be two smooth harmonic almost complex structures, where $M$ is connected. If they agree on an open set $U\subset M$, then they are identical. 
			In particular, a harmonic almost complex structure which is constant on an open set is a constant on the whole manifold $M$.
		\end{theorem}
		\begin{proof}
			We recall Aronszajn's generalization of Carleman's unique continuation theorem \cite[p.248]{Aronszajn-1957}: \emph{Let $A$ be a linear elliptic second-order differential operator defined on a domain $D\subset\R^m$. In $D$ let $u=(u^1,...,u^r)$ be functions satisfying the inequality
				\begin{equation}\label{eq: Aronszajn's ineq}
					\left|Au^\alpha\right|\le C\left(\sum_{i,\beta}\left|\frac{\partial u^\alpha}{\partial x^i}\right|+ \sum_{\beta}|u^\beta|\right).
				\end{equation}
				If $u=0$ in an open set, then $u = 0$ throughout $D$.}
			We apply Aronszajn's inequality \eqref{eq: Aronszajn's ineq} to $(J-J')$ as follows. 
			From \eqref{eq: Harmonic ACS eq} we have 
			$$\begin{aligned}
				\Delta(J-J')&=J\nabla J\cdot\nabla J-J'\nabla J'\cdot\nabla J'\\
				&=J\nabla J\cdot\nabla(J-J')+J\nabla(J-J')\cdot\nabla J'+(J-J')\nabla J'\cdot\nabla J'.
			\end{aligned}$$
			Then, we have
			$$\begin{aligned}
				|\Delta(J-J')|&\le|J\nabla J|\cdot|\nabla(J-J')|+|J\nabla(J-J')|\cdot|\nabla J'|+|J-J'|\cdot|\nabla J'|^2\\
				&\le C\left(|\nabla J|\cdot|\nabla(J-J')|+|\nabla(J-J')|\cdot|\nabla J'|+|J-J'|\cdot|\nabla J'|^2\right).
			\end{aligned}
			$$
			In the open set $U$, slightly shrunk if necessary, the derivatives $\nabla J$ and $\nabla J'$ are bounded. It is easy to see that$$
			|\Delta(J-J')|\le C\left(|\nabla(J-J')|+|J-J'|\right)\qquad \text{in }U.
			$$
			As $J-J'$ vanish in the open set $U$, we have $J-J'\equiv0$ throughout the neighbourhood of $U$. Our conclusion follows from the connectedness of $M$.
		\end{proof}
		
		\begin{proof}[Proof of Proposition \ref{prop: unique continuation}] 
			Note that if $J$ and $J'$ are both smooth, then the claim follows from Theorem \ref{thm: unique continuation}. In the general case, let $\Sigma$ and $\Sigma'$ be the closed set for $J$ and $J'$, respectively, such that $J\in C^\infty(B_1\backslash\Sigma)$, $J'\in C^\infty(B_1\backslash\Sigma')$. Then $J$ and $J'$ are smooth on $B_1\backslash(\Sigma\cup\Sigma')$. By Theorem \ref{thm: ep-regularity} we have $\cal{H}^{m-2}(\Sigma\cup\Sigma')<\infty$ and thus non-disconnecting. Applying Theorem \ref{thm: unique continuation} once again gives $J=J'$ almost everywhere on $B_1\backslash(\Sigma\cup\Sigma')$ and thus on the whole domain.
		\end{proof}
		
		\subsection{Tangent almost complex structure}\label{sec: Tangent almost}
		In this subsection, we shall explore some fundamental properties of tangent almost complex structure. {Here we follow the presentation of He \cite[Section 4.4]{He-2019}. For $a\in M$ such that $B_1(a)\subset M$. Since the exponential map $\exp_a: B_1(0)\subset \R^{2n}\rightarrow B_1(a)\subset M$ is a diffeomorphism. We define the rescaled map $J_{a,r}(y):= J\circ\exp_a(ry)$ and the rescaled metric $g_{a,r}(y):= g\circ\exp_a(ry)$ for $y\in B_{1/r}(0)\subset\R^{2n}$.  A basic observation is that if $J\in W^{1,2}(\cJ_g(B_1(p)))$ is a minimizing harmonic almost complex structure, then the rescaled map $J_{a,r}(y)\in W^{1,2}(\cJ_{g_{a,r}}(B_{1/r}(0)))$ is also a minimizing harmonic almost complex structure.} {For convenience, we still denote $J\circ\exp_a(ry)$ and $g\circ\exp_a(ry)$ by $J(a+ry)$ and $g(a+ry)$, respectively.}
		\begin{definition}\label{def:tangent hacs}
			We say that $T\in W^{1,2}(\mathbb{R}^m,SO(2n)/U(n))$ is a tangent almost complex structure of a minimizing harmonic almost complex structure $J\in W^{1,2}(\cJ_g)$ at the point $a\in B_1$, if there is a sequence $r_j\to 0$ such that $J_{a,r_j}:= J(a+r_j\cdot)\wto T$ in $W^{1,2}$.
		\end{definition}
		
		Based on Theorem \ref{thm: compact minimizing}, we are able to deduce the following basic results for tangent almost complex structures of minimizing harmonic almost complex structures.
		\begin{proposition} \label{prop: tangent map}
			{Let $J\in W^{1,2}_\La(\cJ_g(B_1(p)))$ be a minimizing harmonic almost complex structure for $B_1(p)\subset M$. Suppose there is a sequence $r_i\to 0$ such that
				$$
				J_{i}\wto T\quad\text{weakly in }W^{1,2}(B_{1}(0)),
				$$
				where $J_{i}=J_{a,r_i}:=J(a+r_i\cdot)$ on $a\in B_1(p)$.} Then
			
			{\upshape(i)} $J_{i}\to T$ strongly in $W^{1,2}$. 
			
			{\upshape(ii) (Symmetry of tangent)} $T$ is a minimizing harmonic almost complex structure with $\mathcal{H}^{m-3}(\sing(T))=0$. Moreover, it is $0$-homogeneous with respect to the origin, i.e.
			$$T(x)=T(\lambda x),\quad\text{for all }x\in\R^m \text{ and }\lambda>0.$$
		\end{proposition}
		
		\begin{proof}
			Assertion (i) is a direct consequence of Theorem  \ref{thm: compact minimizing}. Now we prove the homogeneity of $T$ in (ii). 
			
			By assertion (i) and the definition of density function $\Phi_J$, it follows that 
			$$\begin{aligned}
				r^{2-m}\int_{B_r(0)}&|DT|^2(x)dx
				=\lim_{i\to\oo}r^{2-m}\int_{B_r(0)}|D J_i|^2(x)dx\\
				&=\lim_{i\to\oo}(rr_i)^{2-m}\int_{B_{rr_i}(a)}|DJ|^2(y)dy\\
				&=\lim_{t\to0}t^{2-m}\int_{B_{t}(a)}|DJ|^2(y)dy=\Phi_J(a),
			\end{aligned}
			$$which is independent of $r$. Therefore, by applying the  monotonicity formula \eqref{eq: delta-monotonicity formula} to $T$, we obtain 
			$$
			\int_{B_{R}\setminus B_{r}(0)}\rho^{2-m}|\partial_{\rho}J_i|^{2}dx\le R^{2-m}\int_{B_R(0)}|DT|^2-r^{2-m}\int_{B_r(0)}|DT|^2=0.
			$$
			This implies that $T$ is radially invariant in $B_{R}\setminus B_{r}(0)$. Since $R,r$ are arbitrary, this implies that $T$ is radially invariant with respect to the origin. Since $T$ is minimizing, by Theorem \ref{thm: ep-regularity}, $\mathcal{H}^{m-3}(\sing(T))=0$. The proof is complete.
		\end{proof}

		\section{Quantitative $\ep$-regularity theorem}\label{sec:quantitative ep regularity}

		In this section, we develop a stratification theory of the singular set for minimizing harmonic almost complex structures based on the symmetry of tangent harmonic almost complex structures. The original idea is due to Cheeger-Naber \cite{Cheeger-Naber-2013-CPAM} and Naber-Valtorta \cite{Naber-V-2017} in their study of harmonic maps. 
		
		\subsection{Quantitative symmetry and cone splitting principle}

		{We follow the reduction of He \cite[Section 5.1]{He-2019}. Consider a geodesic ball $B_3(p)\subset M$, we identify $B_3(p)$ with an Euclidean ball $B_3(0)$ via the exponential map, $\exp: B_3(0)\rightarrow M$ with $\exp(0)=p$. Again we assume that the injectivity radius of $(M, g)$ is larger than $3$, by a fixed scaling if necessary. We also assume that the metric $g$ is universally close to an Euclidean metric in $B_2(0)$, given by \eqref{eq: g near Euclidean}. Similarly to Section \ref{sec: Tangent almost}, we consider the almost complex structure defined on $B_{\frac{3}{r}}(0)\subset\R^m$ by $J\circ\exp_p(ry)$, still denoted by $J(p+ry)$. 
			
			By the results in Section \ref{sec:partial regularity}, for any sequence $r_i\to0$, $J(p+r_iy)$ converges by subsequence to a tangent almost complex structure $T$ on $\R^m$ which is compatible with the Euclidean metric. Moreover, by Theorem \ref{thm: compact minimizing} and Proposition \ref{prop: tangent map}, $T$ is energy minimizing and homogeneous of degree zero, which gives the classical stratification of the singular set $\sing(J)$ based on homogeneity of its tangent map. For a quantitative version of such stratification, we consider the almost homogeneity property of $J$.}

		\begin{definition}[Symmetry]\label{def: symmetry}
			{Given a measurable map $J\colon \mathbb{R}^m\to N$ we say that}
			
			{
				{\upshape(1)} $J$ is $0$-homogeneous or $0$-symmetric with respect to point $p$ if $J(p+\lambda v)=J(p+v)$ for all $\lambda>0$ and $v\in\mathbb{R}^m$.}
			
			{
				{\upshape(2)} $J$ is $k$-symmetric if it is $0$-homogeneous with respect to the origin,  and is translation invariant with respect to a $k$-dimensional subspace $V\subset\mathbb{R}^m$, i.e.,
				$$J(x+v)=J(x) \qquad \text{for all } x\in\mathbb{R}^m,\, v\in V.$$}
		\end{definition}
		
		{For $J:M\to N$, we say $J$ is $0$-homogeneous with respect to $p$ if $J\circ\exp_p(\lambda v)=J\circ\exp_p(v)$. We say $J$ is $k$-symmetric with respect to $p$ if it is $0$-homogeneous with respect to $p$, and is translation invariant with respect to a $k$-dimensional subspace $V\subset\mathbb{R}^m$, i.e.,
			$$J\circ\exp_x(v)=J(x) \qquad \text{for all } x\in M,\, v\in B_1(0)\cap V.$$}
		
		{If $J\in C^1(\mathbb{R}^m,N)$, then $J$ is $0$-homogeneous at $x=p$ if and only if $\partial_{r_p}J=0$, where $r_p$ is the radial direction centered at $p$; and $J$ is translation invariant with respect to a $k$-dimensional subspace $V\subset\mathbb{R}^m$ if and only if $\partial_v J=0$ for all $v\in V$.}
		
		\begin{definition}[Quantitative  symmetry]\label{def: quantitative symmetry}
			Given a map $J\in L^2(M,N)$, $\ep>0$ and nonnegative integer $k$,  we say that $J$ is $(k,\ep)$-symmetric  on $B_r(x)\subset M$, or simply $B_r(x)\subset M$ is $(k,\ep)$-symmetric, if there exists some $k$-symmetric function $K\colon M\to N$ with respect to $x$ such that
			$$\medint_{B_r(x)}|J(y)-K(y)|^2\leq\ep.$$
		\end{definition}

		Note that $J$ is $(k,\ep)$-symmetric on $B_r(x)$ if and only if the scaled map  $J_{x,r}(y)=J(x+ry)$ is $(k,\ep)$-symmetric  on $B_1(0)$.
		
		According to the formal definition of quantitative symmetry, we may stratify the domain of a function by classifying points according to the degree of symmetry exhibited by the function in infinitesimal neighborhoods of those points. This quantitative stratification partitions the domain into subsets where each subset corresponds to a distinct level of symmetric regularity.
		\begin{definition}[Quantitative stratification]\label{def: quantitative stratification}
			For any map $J\in L^2(M,N)$, $r,\ep>0$ and $k\in\{0,1,\cdots,m\}$, we define the $k$-th quantitative singular stratum $S_{\ep,r}^k(J)\subset M$ as 
			$$S^k_{\ep,r}(J)\equiv\left\{x\in M: J \text{ is not } (k+1,\ep)\text{-symmetric on }B_s(x) \text{ for any } r\le s<1\right\}.$$
			Furthermore, we set
			$$
			S^k_\ep(J):=\bigcap_{r>0}S^k_{\ep,r}(J)
			\quad\text{ and }\quad 
			S^k(J)=\bigcup_{\ep>0}S^k_\epsilon(J).
			$$
		\end{definition}
		
		Note that,
		$$k'\leq k\text{ or }\ep'\geq\ep\text{ or } r'\leq r\ \Longrightarrow \ S^{k'}_{\ep',r'}(J)\subseteq S^k_{\ep,r}(J).$$
		In particular, we have
		$$S^0(J)\subset S^1(J)\subset\cdots\subset S^m(J)=M.$$
		
		Now we apply this method to singular sets of minimizing harmonic almost complex structures.
		\begin{lemma}\label{lemma: stratification for SHACS} Suppose $J\in W^{1,2}(\cJ_g(B_1(p)))$ is a minimizing harmonic almost complex structure on $B_1(p)\subset M$, then
			$$S^k(J)=\{x\in B_1(p): \text{no tangent almost complex structure of } J \text{ at } x \mbox{ is } (k+1)\mbox{-symmetric}\}.$$
			Consequently, we have
			$$S^0(J)\subset S^1(J)\subset\cdots\subset S^{m-3}(J)\subset S^{m-2}(J)\subset S^{m-1}(J)=\sing(J).$$
		\end{lemma}

		\begin{proof} 
			First of all, we let
			$$\Sigma^k(J)=\{x\in B_1(p): \text{no tangent almost complex structure of }J \text{ at } x \mbox{ is } (k+1)\mbox{-symmetric}\}$$at this moment. 		
			Suppose $y\in S^k(J)$, then $y\in S^k_\ep(J)$ for some $\ep>0$. Thus, for any $(k+1)$-symmetric map $K\in L^2$ and $r>0$, we have
			$$\medint_{B_1(0)}|J_{y,r}-K|^2dx\geq\ep.$$
			If $T$ is a tangent map of $J$ at $y$, then there exists a sequence $r_i\to 0$ such that $J_{y,r_i}\to T$ strongly in $W^{1,2}_{\loc}$. This gives $$\medint_{B_1(0)}|K-T|^2dx\geq\ep,$$
			which implies that $T$ is not $(k+1)$-symmetric. This implies  $S^k(J)\subset \Sigma^k(J)$.
			
			Conversely, suppose $y\notin S^k(J)$. Then there exist sequences $r_i>0$ and $(k+1)$-symmetric maps $K_i$ from $\R^m$ to $ SO(2n)/U(n)$ such that
			$$\medint_{B_1(0)}|J_{y,r_i}-K_i|^2dx\leq i^{-1}.$$
			We can select subsequences $\{J_{x,r_i}\} \rightharpoonup T$ in $W^{1,2}$ and $K_i\rightharpoonup K$ in $L^2(B_1(0))$. Then  by the weak lower semi-continuity of $L^2$-norm, we obtain
			$$
			\int_{B_1(0)}|T-K|^2dx\leq\liminf_{i\to\oo}\int_{B_1(0)}|J_{y,r_i}-K_i|^2dx=0.
			$$
			Moreover, by the compactness of symmetric maps, $K$ is $(k+1)$-symmetric.
			
			If $r_i\to0$,  then $T$ is a tangent map and thus is $(k+1)$-symmetric, which shows that $x\not\in \Sigma^k(J)$. 
			
			If $\lim\limits_{i\to\oo}r_i=r>0$, we claim that $J(z)=K(z-y)$ for almost every $z\in B_\delta(y)$, where $\delta\in(0,r]$ is a positive radius. Indeed, we notice that
			$$
			\medint_{B_{r_i}(y)}|J(z)-K_i(z-y)|^2dz=\medint_{B_1(0)}|J_{y,r_i}(z)-K_i(z)|^2dz\leq i^{-1},
			$$
			which implies that there exists $\delta\in(0,r]$ such that
			$$
			\medint_{B_\delta(y)}|J(z)-K(z-y)|^2dz\le\liminf_{i\to\infty}\medint_{B_1(0)}|J_{y,r_i}(z)-K_i(z)|^2dz=0,
			$$
			using the lower semi-continuity of $L^2$-norm. Thus, all tangent maps of $J$ at $y$ are $(k+1)$-symmetric and thus $x\notin \Sigma^k(J)$. The proof is complete.
		\end{proof}

		Now we are ready to establish quantitative symmetry and cone splitting principle for minimizing harmonic almost complex structures $J$. Recall from \eqref{eq: mono energy} that 
		$$
		\Phi_{J}(a, r)=e^{C_0r}r^{2-m}\int_{B_r(a)}|DJ|^2+C_0r
		$$
		
		\begin{proposition}\label{prop: quantitative symmetry}
			{For any $\ep>0$, there exists $\delta_0=\delta_0(m,g,\La,\ep)$ such that, if $J\in W^{1,2}_\La(\cJ_g(B_3(p)))$ is a minimizing harmonic almost complex structure on $B_3(p)\subset M$} with
			$$\Phi_{J}(x,r)-\Phi_{J}(x,r/2)<\delta_0$$
			for some $x\in B_1(p)$ and $0<r<1$, then $J$ is $(0,\ep)$-symmetric on $B_r(p)$.
		\end{proposition}
		\begin{proof}
			We prove by contradiction. Suppose for some $\ep>0$ there is a sequence of  minimizing harmonic almost complex structure {$J\in W^{1,2}_\La(\cJ_g(B_3(p)))$} satisfying
			$$\Phi_{J_i}(x_i,r_i)-\Phi_{J_i}(x_i,r_i/2)<i^{-1},$$
			but is not $(0,\ep)$-symmetric on $B_{r_i}(x_i)$.
			
			{Let $\bar{J}_i(z)=J_i(x_i+r_iz)(=J_i\circ\exp_{x_i}(r_iz))$.} 
			Since $\{\bar{J}_i\}$ is uniformly bounded in $W^{1,2}$, we can select a subsequence, still denote by $\{\bar{J}_i\}$, such that  $\bar{J}_i\rightharpoonup J$ weakly in $W^{1,2}$ and strongly in $L^2$. Then, by Theorems \ref{thm: ep-regularity} and \ref{thm: compact minimizing}, $J$ is minimizing in $B_1$ with $\dim_{\mathcal{H}}\sing(J)\le m-3$.
			By the monotonicity formula \eqref{eq: monotonicity formula}, for all $i\in \N$, we have
			$$
			\begin{aligned}
				\int_{B_1\backslash B_{1/2}}{|y\cdot D\bar{J}_i|^2}dy =&
				\int_{B_{r_i}\backslash B_{\frac{r_i}{2}}(x_i)}|(x-x_i)\cdot DJ_i|^2r_i^{-m}dx\\
				\le& C\int_{B_{r_i}\backslash B_{\frac{r_i}{2}}(x_i)}|x-x_i|^{2-m}|\partial_\rho J_i|^2dx\\
				\le& C(\Phi_{J_i}(x_i,r_i)-\Phi_{J_i}(x_i,r_i/2))<Ci^{-1}.
			\end{aligned}
			$$
			Sending $i\to \infty$ and using the strong convergence of $\bar{J}_i\rightarrow J$ in $L^2$, we conclude that $J$ is radially invariant on $B_1\backslash B_{1/2}(0)$,  and thus is homogeneous on $\R^m$ by the unique continuation property (see proposition \ref{prop: unique continuation}).
			
			Finally, since $\bar{J}_i\rightarrow J$ strongly in $L^2$, we have$$
			\medint_{B_{r_i}(x_i)} |J_i(y)-J((y-x_i)/r_i)|^2dy=\medint_{B_{1}}|\bar{J}_i-J|^2\to 0\qquad \text{as }i\to \infty, 
			$$
			contradicting with the assumption that $J_i$ is not $(0,\ep)$-symmetric on $B_{r_i}(x_i)$. The proof is thus complete. 
		\end{proof}
		
		In order to extend proposition \ref{prop: quantitative symmetry} to higher order symmetric cases, we introduce the notion of quantitative frame for affine subspaces as in \cite{Naber-V-2018}.
		
		\begin{definition}
			{Let $\{y_i\}_{i=0}^k \subset B_1(p)\subset M$ and $\rho>0$. We say that these points $\rho$-effectively span a $k$-dimensional affine subspace if for all $i=1,\cdots,k$,
				$$\dist(z_i,z_0+\Span\{z_1-z_0,\cdots,z_{i-1}-z_0\})\geq2\rho,$$
				where $z_j=\exp^{-1}(y_j)\in\R^m$ for $j=0,...,k$.
				More generally, a set $F\subset B_1(p)$ is said to $\rho$-effectively span a $k$-dimensional affine subspace, if there exist points $\{y_i\}_{i=0}^k\subset F$ which $\rho$-effectively spans a $k$-dimensional affine subspace.}
		\end{definition}
		
		To establish the quantitative version of cone splitting principle, we recall the observation in \cite{Cheeger-Naber-2013-CPAM}.
		
		\begin{proposition}[Cone splitting principle]\label{prop: cone splitting}
			{If $J: \mathbb{R}^{\ell_1}\rightarrow\mathbb{R}^{\ell_2}$ is $k$-symmetric with respect to a $k$-plane $V^k,$ and is $0$-symmetric at a point $z\not\in V^k$, then $J$ is $(k+1)$-symmetric with respect to  the $(k+1)$-plane $V^{k+1}:=\Span\{z,V^k\}.$} 
		\end{proposition}
		
		The following result extends Proposition \ref{prop: quantitative symmetry} to higher order quantitative symmetry.
		
		\begin{proposition}[Quantitative cone splitting]\label{prop: quant cone splitting}
			{For any $\ep,\rho>0$, there exists $\delta_1=\delta_1(m,g,\ep,\rho)>0$ such that, if $J\in W^{1,2}_\La(\cJ_g(B_3(p)))$ is a minimizing harmonic almost complex structure on $B_3(p)\subset M$, and for some $\{x_j\}_{j=0}^k\subset B_1(p)$} there holds
			
			{\upshape(a)} $\Phi_{J}(x_j,r)-\Phi_{J}(x_j,r/2)<\delta_1 $ for some $0<r<1$ and for all $j=0,\cdots,k$,
			
			{\upshape(b)} $\{x_j\}$ $\rho$-effectively spans a $k$-dimensional affine plane $V$,
			
			\noindent then $J$ is $(k,\ep)$-symmetric on $B_1(p)$.
		\end{proposition}

		\begin{proof}
			Suppose by contradiction that there is a sequence of minimizing harmonic almost complex structure {$J^i\in W^{1,2}_\La(\cJ_g(B_3(p)))$}, $\ep_0,\rho_0>0$, $\{x_0^i,\cdots,x_k^i\}$ and a sequence of $k$-dimensional affine plane $V^i$ such that $\{x_0^i,\cdots,x_k^i\}$  $\rho_0$-effectively spans $V^i$ and $\Phi_{J^i}(x^i_j,r)-\Phi_{J^i}(x^i_j,r/2)<\delta_i$, where $\delta_i\to 0$ is chosen according to proposition \ref{prop: quantitative symmetry}, so that $J^i$ is $(0,i^{-1})$-symmetric at $B_{i^{-1}}(x_j^i)$ for each $j$. But $J^i$ is not $(k,\ep_0)$-symmetric in $B_1(0)$ for all $i\geq1.$ 
			
			Assume $J^i\rightharpoonup J$ weakly in $W^{1,2}(B_1(0))$ and strongly in $L^2(B_1(0))$, $x_j^i\rightarrow x_j\in B_1(0)$ and $V^i\ri V$ as $i\rightarrow\infty.$
			Applying the monotonicity formula \eqref{eq: monotonicity formula} to each $J^i$, we obtain
			$$\int_{B_r(x^i_j)\backslash B_{r/2}(x_j^i)}|(x-x^i_j)\cdot DJ^i|^2\leq C\left(\Phi_{J^i}(x^i_j,r)-\Phi_{J^i}(x^i_j,r/2)\right)<C\delta_i.$$
			Sending $i\ri\oo$, we get
			$$\int_{B_r(x_j)\backslash B_{r/2}(x_j)}|(x-x_j)\cdot DJ|^2=0.$$
			This shows that $J$ is radially symmetric on ${B_r(x_j)\backslash B_{r/2}(x_j)}$ for each $j\in \{0,\cdots,k\}$. The unique continuation property (see proposition \ref{prop: unique continuation}) and the cone splitting principle (see proposition \ref{prop: cone splitting}) imply $J$ is symmetric with respect to the $k$-dimensional affine plane spanned by $\{x_j\}_{j=0}^k$. This contradicts with the assumption that $J^i$ is not $(k,\ep_0)$-symmetric in $B_1(0)$ since $J^i\to J$ strongly in $L^2$.
		\end{proof}

		\subsection{Basic properties of quantitative stratum} 
		
		The following proposition shows that under some conditions, the singular set satisfies the  one-side Reifenberg approximating property.
		\begin{proposition}\label{prop: stratum property-1} 
			For any $\ep,\rho>0$, there exists $\delta_2=\delta_2(m,g,\Lambda,\ep,\rho)>0$ such that: If $J\in W^{1,2}_\La(\cJ_g(B_3(p)))$ is a minimizing harmonic almost complex structure on $B_3(p)\subset M$ and the set $$
			F:=\{y\in  B_2(p): \Phi_J(y,2)-\Phi_J(y,\rho)<\delta_2\}$$
			$\rho$-effectively spans a $k$-dimensional affine plane $V$, then $$S^k_{\ep,\delta_2}(J)\cap B_1(p)\subset \exp_p(B_{2\rho}(V)).$$
		\end{proposition}

		We need the following lemma to prove proposition \ref{prop: stratum property-1}.  This result shows the relationship between almost translation invariant property and  almost $(k+1)$-symmetric property. Here and after, $|P\cdot Du|^2=\sum\limits_{i=1}^{k+1}|D_{e_i}u|^2$ with $\{e_i\}_{i=1}^{k+1}$ being an  orthonormal basis of a $(k+1)$-dimensional subspace $P$.
		\begin{lemma}\label{lemma: almost translation invariant}
			For any $\ep>0$, there exists $\delta_3=\delta_3(m,g,\Lambda,\ep)>0$ such that if $J\in W^{1,2}_\La(\cJ_g(B_3(p)))$ is a minimizing harmonic almost complex structure on $B_3(p)\subset M$ satisfying
			\begin{align}\label{eq: almost translation invariant}
				\int_{B_1(p)}|P\cdot D J|^2<\delta_3
			\end{align}
			for some $(k+1)$-dimensional subspace $P$, then 
			$S^k_{\ep,\overline{r}}(J)\cap B_{1/2}(p)=\emptyset$ for  $\overline{r}=\delta_3^{\frac{1}{2(m-2)}}$. In particular, $0\not\in \exp^{-1}S^k_{\ep,\overline{r}}(J)$, i.e. $p\not\in S^k_{\ep,\overline{r}}(J)$.
		\end{lemma}
		The proof of Lemma \ref{lemma: almost translation invariant} is similar to Naber-Valtorta \cite[Lemma 4.8]{Naber-V-2018}. 
		\begin{proof}	
			Identify $B_3(p)\subset M$ with $B_3(0)\subset \R^m$ via the exponential map $\exp$, and still denote the pullback maps $\exp_p^*J$ and $\exp_p^*g$ by $J$ and $g$, respectively.
			We claim that there is a constant $C(m,\Lambda)>0$ such that for every $x\in B_{1/2}(p)$, there exists $r_x\in[\overline{r},1/2]$ such that
			\begin{equation}\label{eq: lemma 3-11-1}
				\Phi_J(x,r_x)-\Phi_J(x,r_x/2)<\frac{C(m,\Lambda)}{-\log\delta_3}.
			\end{equation}
			Indeed, if this is not true, then for $J\in W^{1,2}_\La(\cJ_g(B_3(p)))$ we may assume by choosing a good radius that $\Phi_J(x,1/2)\le C_0(m,\Lambda)$. Then 
			\begin{align*}
				C_0(m,\Lambda)\geq\Phi_J(x,1/2)&\geq\sum_{i=1}^{-\log\overline{r}+1}(\Phi_J(x,2^{-i})-\Phi_J(x,2^{-i-1}))\geq c(m)C(m,\Lambda),
			\end{align*}
			which is impossible if we let $C(m,\Lambda)=\frac{2C_0(m,\Lambda)}{c(m)}$. This proves the claim. 
			
			Suppose by contradiction that, there exist $\ep>0$ and a sequence of minimizing harmonic almost complex structure $J\in W^{1,2}_\La(\cJ_g(B_3))$ with subspaces $P_i$, $\de_{3,i}\to 0$, $x_i\in B_{1/2}$ and $r_{x_i}\in[\bar{r}_i,1]$, such that $J_i$ is not $(k+1,\ep)$-symmetric on $B_{r_{x_i}}(x_i)$, \eqref{eq: lemma 3-11-1} holds for $x=x_i$, where $\bar{r}_i=\de_{3,i}^{1/2(m-2)}$. Note also that by the definition of $r_{x_i}$, we have
			\begin{align}\label{eq: lemma 3-11-2}
				r_{x_i}^{2-m}\int_{B_{r_{x_i}}(x_i)}|P_i\cdot Du_i|^2<r_{x_i}^{2-m}\delta_{3,i}\leq\de_{3,i}^{1/2}.
			\end{align} 
			After a simple rotation, we may assume that the $(k+1)$-dimensional subspaces $P_i$ are equal to $P$ for all $i$. 
			
			Let $T_i(x)=J_i(x_i+r_{x_i}x)$. Then we may assume that $T_i$ converges weakly in $W^{1,2}$ and strongly in $L^2$ to $T$. By theorem \ref{thm: compact minimizing}, $T$ is a minimizing harmonic almost complex structure and $T_i\to T$ strongly in $W^{1,2}$. It follows that $T$ is $0$-symmetric by unique continuation property (see proposition \ref{prop: unique continuation}) and invariant with respect to the $(k+1)$-dimensional subspace $P$ by \eqref{eq: lemma 3-11-2}: 
			$$
			0=\lim_{i\to \infty}r_{x_i}^{2-m}\int_{B_{r_{x_i}}(x_i)}|P\cdot D J_i|^2=\lim_{i\to \infty}\int_{B_1}|P\cdot DT_i|^2=\int_{B_1}|P\cdot DT|^2.
			$$
			This implies that $T_i$ is $(k+1,\ep)$-symmetric on $B_{1}$ for $i$ big enough, or equivalently, $J_i$ is $(k+1,\ep)$-symmetric with $B_{r_i}(x_i)$, which is a contradiction.
		\end{proof}
		
		Now we can prove Proposition \ref{prop: stratum property-1}.
		\begin{proof}[Proof of Proposition \ref{prop: stratum property-1}]
			Suppose $\{y_{j}\}_{j=0}^{k}\subset F$ is a $\rho$-independent frame
			that spans the affine subspace $V$, i.e., 
			$$
			V=\exp_p^{-1}y_{0}+{\rm span}\{\exp_p^{-1}y_{i}-\exp_p^{-1}y_{0}\}_{i=1}^{k}.
			$$
			Let {$x_{0}\in B_{1}(p)\backslash \exp_pB_{2\rho}(V)$} and $\de_{2}>0$ to be
			determined later. Then we prove that $x_{0}\not\in S_{\ep,\de_{2}}^{k}(J)\cap B_{1}(p)$. 
			By the definition of $F$, there holds $\Phi_{J}(y_{i},2)-\Phi_{J}(y_{i},\rho)<\de_{2}$
			for all $0\le i\le k$. We can select a positive $r$ small enough such that  $B_{r}(x_{0})\subset B_{2}(y_{i})\backslash B_{\rho}(y_{i})$
			for every $i$. Thus, by the monotonicity formula \eqref{eq: monotonicity formula}, for each $0\le i\le k$, it holds
			$$
			\int_{B_{r}(x_{0})}|(z-y_{i})\cdot DJ|^{2}dz\le\int_{B_{2}(y_{i})\backslash B_{\rho}(y_{i})}|(z-y_{i})\cdot DJ|^{2}dz\le C\de_{2}.
			$$
			Consequently, we obtain
			$$
			\begin{aligned}
				&\int_{B_{r}(x_{0})}|(y_{i}-y_{0})\cdot DJ|^{2}\\
				\le&2\int_{B_{2}(y_{i})\backslash B_{\rho}(y_{i})}|(z-y_{i})\cdot DJ|^{2}+2\int_{B_{2}(y_{i})\backslash B_{\rho}(y_{i})}|(z-y_{0})\cdot DJ|^{2}\\
				\le&2C\de_{2}.
			\end{aligned}
			$$
			Since $\{y_{j}\}_{j=0}^{k}$ is $\rho$-independent, we conclude that
			\begin{equation}\label{eq: hatV nabla J small}
				\int_{B_{r}(x_{0})}|\hat{V}\cdot DJ|^{2}\le C(m,\rho)\de_{2},
			\end{equation}
			where {$\hat{V}={\rm span}\{\exp_p^{-1}y_{i}-\exp_p^{-1}y_{0}\}_{i=1}^{k}$}.
			On the other hand, for each $z\in B_{r}(x_{0})\subset B_{2}\backslash B_{\rho}(V)$,
			let {$
				\pi_{V}(z)=\exp_p\left(\exp_p^{-1}y_{0}+\sum_{i=1}^{k}\al_{i}(z)(\exp_p^{-1}y_{i}-\exp_p^{-1}y_{0})\right)
				$} be the orthogonal projection of $z$ in $V$. Then $|z-\pi_{V}(z)|\ge\rho$,
			$|\al_{i}(z)|\le C(m,\rho)$ and
			$$
			\int_{B_{r}(x_{0})}|(z-\pi_{V}(z))\cdot DJ|^{2}dz\le C\sum_{i=0}^{k}\int_{B_{2}(y_{i})\backslash B_{\rho}(y_{i})}|(z-y_{i})\cdot DJ|^{2}\le C(m,\rho)\de_{2}.
			$$
			Thus, by setting $h(z)=(z-\pi_{V}(z))/|z-\pi_{V}(z)|$, it follows
			that
			$$
			\begin{aligned}\int_{B_{r}(x_{0})}|h(x_{0})\cdot DJ|^{2} & \le2\int_{B_{r}(x_{0})}|h(z)\cdot DJ|^{2}+2\int_{B_{r}(x_{0})}|(h(z)-h(x_{0}))\cdot DJ|^{2}\\
				& \le C\de_{2}+Cr^{2}\int_{B_{r}(x_{0})}|DJ|^{2}\\
				& \le C\de_{2}+C\Lambda r^{m}.
			\end{aligned}
			$$
			Now we choose $r=r(m,\La,\rho)\ll\rho$ such that $C\Lambda r^{m}\le C\de_{2}$.
			This gives
			$$
			\int_{B_{r}(x_{0})}|h(x_{0})\cdot DJ|^{2}\le2C\de_{2}.
			$$
			Together with \eqref{eq: hatV nabla J small}, for the $(k+1)$-dimensional subspace $P:=\hat{V}\oplus \R h(x_0)$, we find that
			$$
			\int_{B_{r}(x_{0})}|P\cdot DJ|^{2}\le C(m,\rho,\La)\de_{2}.
			$$
			By choosing $\de_{2}=\de_{2}(m,\rho,\La)$ sufficiently small so that we can apply Lemma \ref{lemma: almost translation invariant} to conclude that $J$ is $(k+1,\ep)$-symmetric in $B_{r}(x_{0})$. This implies {$S^k_{\ep,\delta_2}(J)\cap B_1\subset \exp_pB_{2\rho}(V)$}.
		\end{proof}

		The next result shows that $\Phi_J$ remains almost constant on all pinched points.
		\begin{lemma}\label{lemma: Phi_J almost const}
			Let $\rho,\eta>0$ be fixed and {$J\in W^{1,2}_\La(\cJ_g(B_3(p)))$} is a minimizing harmonic almost complex structure. Let  
			{$$E:=\sup_{y\in B_2(p)}\Phi_J(y,1).$$ }
			Then there exists $\de_4=\de_4(m,g,\Lambda,\rho,\eta)>0$ such that, if the set 
			$$F:=\{y\in B_2(p):\Phi_J(y,\rho)>E-\delta_4\}$$
			$\rho$-effectively spans a $k$-dimensional affine subspace $L\subset\mathbb{R}^m$, then
			{$$\Phi_J(x,\rho)\ge E-\eta \qquad \text{for all }  x\in \exp_p(L\cap B_2(0)).$$}
		\end{lemma}
		\begin{proof}
			We prove by contradiction. Suppose {$J\in W^{1,2}_\La(\cJ_g(B_3))$} is a sequence of minimizing harmonic almost complex structure with
			$$\sup_{y\in B_2(p)}\Phi_{J_i}(y,1)\leq E.$$
			{For each $i\ge 1$, we assume that there exists $\{y^i_j\}_{j=0}^k\subset F_i:=\{y\in B_2(p):\Phi_{J_i}(y,\rho)>E-i^{-1}\}$  spanning $\rho$-effectively a $k$-dimensional affine subspace $L_i\subset\mathbb{R}^m$, and there exists $x_i\in \exp_p(L_i\cap B_2(0))$} such that
			\begin{equation}\label{eq: Phi_J almost const-1}
				\Phi_{J_i}(x_i,\rho)\leq E-\eta.
			\end{equation}
			It follows from the assumption that
			\begin{equation}\label{eq: Phi_J almost const-2}
				\Phi_{J_i}(y^i_j,1)-\Phi_{J_i}(y^i_j,\rho)<1/i,\qquad \text{for all }\,i\ge 1,\, 0\le j\le k. 
			\end{equation}
			Without loss of generality, we assume that for each $0\le j\le k$, 
			$$y^i_j\to y_j,\quad x_i\to x\quad\text{and}\quad L_i\to L$$
			as $i\to \infty$ with {$x\in \exp_pL$} and $J_i\wto J$ in $W^{1,2}$.
			
			By Theorem \ref{thm: compact minimizing}, we have $J_i\to J$ strongly in $W^{1,2}$ and $J$ is a minimizing harmonic almost complex structure. By Proposition \ref{thm: monotonicity formula}, the function $r\mapsto \Phi_{J}(y,r)$ is monotonically nondecreasing. Sending $i\to \infty$ in \eqref{eq: Phi_J almost const-1} and \eqref{eq: Phi_J almost const-2}, we obtain
			$$
			\Phi_{J}(y_j,1)=\Phi_{J}(y_j,\rho)=E,
			\quad \text{for all }\, 0\le j\le k
			$$
			and
			\begin{equation}\label{eq: loss of energy}
				\Phi_J(x,\rho)\le E-\eta.
			\end{equation}
			Moreover, we know that $J$ is translation invariant along $L$. Hence {$\Phi_J|_{\exp_p(L\cap B_2(0))}\equiv \Phi_J(y_j,1)=E$}, which contradicts with \eqref{eq: loss of energy}. The proof is complete.
		\end{proof}
		
		The following technical result shows the almost symmetry under certain pinching condition.
		\begin{lemma}\label{lemma: pinching almost symmetry}
			For any $\ep,\rho>0$, there exists $\de_5=\de_5(m,g,\Lambda,\rho,\ep)>0$ satisfying: if {$J\in W^{1,2}_\La(\cJ_g(B_3(p)))$ is a minimizing harmonic almost complex structure on $B_3(p)\subset M$ and $\Phi_J(p,1)-\Phi_J(p,1/2)<\de_5$}, and there is a point $y\in B_3(p)$ such that
			
			{\upshape (1)}. $\Phi_J(y,1)-\Phi_J(y,1/2)<\de_5$,
			
			{\upshape (2)}.  $J$ is not $(k+1,\ep)$-symmetric on $B_r(y)$ for some $r\in[\rho,2]$,\\
			then $J$ is not $(k+1,\ep/2)$-symmetric on $B_r(p)$. 
		\end{lemma}
		In particular, under condition {\upshape (1)} of Lemma \ref{lemma: pinching almost symmetry}, {$y\in S^k_{\ep,\rho}(J)\cap B_3(p) \Rightarrow p\in S^k_{\ep/2,\rho}(J)$}.

		\begin{proof}
			Suppose by contradiction that there is a sequence {$\{J_i\}\subset W^{1,2}_\La(\cJ_g(B_3(p)))$} of minimizing harmonic almost complex structure satisfying {$\Phi_{J_i}(p,1)-\Phi_{J_i}(p,1/2)\leq i^{-1},$} and there exists a sequence $\{y_i\}\in B_3(p)$ such that $\Phi_{J_i}(y_i,1)-\Phi_{J_i}(y_i,1/2)\leq i^{-1}$ and that for each $i\in \N$, $J_i$ is not $(k+1,\ep)$-symmetric on $B_r(y_i)$, but it is $(k+1,\ep/2)$-symmetric on $B_r(0)$. 
			Then, there exists a sequence of $(k+1)$-symmetric maps {$K_i\in\Gamma(B_3(p), TM\otimes T^*M)$ such that for any $i\in\N$,
				$$\medint_{B_r(p)}|J_i-K_i|^2\leq\ep/2.$$}
			
			Up to a subsequence if necessary, we may assume that $y_i\ri y\in \overline{B}_3(p)$, $J_i\rii J$ in $W^{1,2}$, and $K_i\rii K$ in $L^2$. By Theorem \ref{thm: compact minimizing} and the unique continuation property Proposition \ref{prop: unique continuation}, we see that $J_i\to J$ strongly in $W^{1,2}$, where $J$ is a minimizing harmonic almost complex structure and is homogeneous with respect to $p$ and $y$. By the compactness of symmetric functions, $K$ is also $(k+1)$-symmetric. Thus, we have
			{$$\medint_{B_r(p)}|J-K|^2\leq\limsup_{i\ri\oo}\medint_{B_r(p)}|J-K_i|^2\leq\ep/2.$$}
			If $y=0$, we get
			\begin{align*}
				&\lim_{i\ri\oo}\medint_{B_r(y_i)}|J_i(x)-K(x-y_i)|^2\\
				\leq&\lim_{i\ri\oo}\medint_{B_r(y_i)}|J_i(x)-K(x)|^2dx+\lim_{i\ri\oo}\medint_{B_r(y_i)}|K(x)-K(x-y_i)|^2dx\\
				\leq&\ep/2.
			\end{align*}
			If $y\neq0,$ then by the invariance of $K$ with respect to $\R y$, we get
			\begin{align*}
				&\lim_{i\ri\oo}\medint_{B_r(y_i)}|J_i(x)-K(x-y_i)|^2=\lim_{i\ri\oo}\medint_{B_r(0)}|J_i(y_i+x)-K(x)|^2dx\\
				\leq& 2\lim_{i\ri\oo}\medint_{B_r(0)}|J_i(y_i+x)-J(y+x)|^2dx+2\medint_{B_r(0)}|J(y+x)-K(x)|^2dx\\
				\leq& \ep.
			\end{align*}
			In both cases, we arrive at a contradiction.	
			\end{proof}
			
			\begin{remark}\label{rmk:on Lemma 313}
				In Section \ref{sec:covering lemma}, we shall repeatedly use the following variant of Lemma \ref{lemma: pinching almost symmetry}: Suppose {$J\in W^{1,2}_\La(\cJ_g(B_3(p)))$ is a minimizing harmonic almost complex structure on $B_3(p)\subset M$} and $\Phi_J(p,1)-\Phi_J(p,1/2)<\de_5$. If there is some $y\in B_3(p)$ with $\Phi_J(y,1)-\Phi_J(y,1/2)<\de_5$, then $y\in S^k_{\ep,\rho}(J)\cap B_3(p)\Rightarrow p\in S^k_{\ep/2,\rho}(J)$.
			\end{remark}

			\section{Reifenberg theorems and estimates of Jones' number}\label{sec:reifenberg}
			In this section, we recall the Reifenberg type results obtained by Naber-Valtorta \cite{Naber-V-2017}. 
			
			We first recall the definition of Jones' number $\beta_2$, which quantifies how close the support of a measure $\mu$ is to a $k$-dimensional affine subspace.
			
			\begin{definition}\label{def: Jones number}
				Let $\mu$ be a nonnegative Radon measure on $B_3$. Fix $k\in\N$. The $k$-dimensional Jones' $\beta_2$ number is defined as
				\begin{align*}
					\beta^k_{2,\mu}(x,r)^2:=\inf_{V\su\R^m}\int_{B_r(x)}\frac{d^2(y,V)}{r^2}\frac{d\mu(y)}{r^k},\qquad \text{when }\,B_r(x)\subset B_3,
				\end{align*}
				where the infimum is taken over  all  $k$-dimensional affine subspaces $V$.
			\end{definition}
			
			We state two versions of the quantitative Reifenberg theorems from \cite{Naber-V-2017}.
			\begin{theorem}[{Discrete-Reifenberg, \cite[Theorem 3.4]{Naber-V-2017}}]\label{thm: discrete-Reifenberg}
				There exist two constants $\de_6=\de_6(m)$ and $C_R(m)$ such that the following property holds. Let $\{B_{r_x}(x)\}_{x\in \mathcal{C}}\su B_2(0)\su\R^m$ be a family of pairwise disjoint balls with centers in $\mathcal{C}\subset B_1(0)$ and let $\mu\equiv\sum_{x\in \mathcal{C}}\omega_kr_x^k\de_x$ be the associated measure. If for every ball $B_r(x)\su B_2,$ there holds
				\begin{align}\label{eq: multiscale appro-1}
					\int_{B_r(x)}\left(\int_0^r\beta^k_{2,\mu}(y,s)^2\frac{ds}{s}\right)d\mu(y)<\de_6^2r^k,
				\end{align}
				then we have the uniform estimate
				\begin{align*}
					\sum_{x\in \mathcal{C}}r^k_x< C_R(m).
				\end{align*}
			\end{theorem}

			\begin{theorem} [{Rectifiable-Reifenberg, \cite[Theorem 3.3]{Naber-V-2017}}]\label{thm: rectifiable-Reifenberg}
				There exist constants $\de_7=\de_7(m)$ and $C=C(m)$ satisfying the following property. Assume that $S\subset B_2\su\R^m$ is $\HH^k$-measurable, and for each $B_r(x)\subset B_2$ there holds
				\begin{align}\label{eq: multiscale appro-2}
					\int_{S\cap B_r(x)}\left(\int_0^r\beta^k_{2,\HH^k|_S}(y,s)^2\frac{ds}{s}\right)d\HH^k(y)<\de_7^2r^k.
				\end{align}
				Then $S\cap B_1$ is  $k$-rectifiable, and   $\HH^k(S\cap B_r(x))\leq C r^k$ for each $x\in S\cap B_1$.
			\end{theorem}
			
			The two conditions \eqref{eq: multiscale appro-1} and \eqref{eq: multiscale appro-2} are usually called multiscale approximation conditions. The key problem in application is how to establish an estimate for the number $\beta_2$.  Similar to the case of harmonic maps \cite[Theorem 7.1]{Naber-V-2017},  we need an $L^2$ subspace approximation theorem for almost complex structure.
			
			For $x\in B_1(p)$ and $r>0,$ we denote for brevity
			$$
			W_r(x):= W_{r,10r}(x)=\int_{B_{10r}(x)\ba B_r(x)}\frac{|(y-x)\cdot DJ(y)|^2}{|y-x|^m}dy\geq0.
			$$
			for an almost complex structure {$J\in W^{1,2}(\cJ_g(M))$}. 
			
			The following result can be proved by the argument of Naber-Vatorta \cite{Naber-V-2017}, whose proof we shall represent in the appendix. 
			\begin{theorem}\label{thm: L2 subspace approximation}
				Fix $\ep>0,0<r\leq1$ and $x\in B_1(p).$ There exist $C(m,g,\Lambda,\ep)>0$ and $\de_8>0,$ such that, if {$J\in W^{1,2}_\La(\cJ_g(B_{12}(p)))$} is a minimizing harmonic almost complex structure defined, and is $(0,\de_8)$-symmetric on $B_{10r}(x)$ but not $(k+1,\ep)$-symmetric, then for any nonnegative finite measure $\mu$ on {$B_r(0)=\exp_p^{-1}(B_r(p))$} we have
				\begin{equation}\label{eq: L2 subspace approximation}
					\beta^k_{2,\mu}(x,r)^2\leq Cr^{-k}\int_{B_r(x)}W_r(y)d\mu(y).
				\end{equation}
			\end{theorem}

			\section{Covering lemma}\label{sec:covering lemma}
			In this section we establish a covering lemma which aims to give an estimate of the energy at singular points. The original idea is from \cite[Section 6.2]{Naber-V-2018} and here we shall follow the slightly modified presentation in \cite[Section 5]{GJXZ-2024}. 
			
			Our main result of this section is to prove the following covering lemma. 
			
			\begin{lemma}[Main covering Lemma]\label{lemma: main covering lemma}
				{Let $J\in W^{1,2}_\La(\cJ_g(B_3(p)))$ be a minimizing harmonic almost complex structure on $B_3(p)\subset M$.} For any $\ep>0$ and $0<r<R\leq1$, there exist constants $\delta=\de(m,\La,\ep)>0$ and $C(m)$ with the following property. 
				
				For any ${S}\su S^k_{\ep,\de r},$ there exists a finite covering of $S\cap B_R(p)$ such that
				\begin{equation}\label{eq: main cover property}
					S\cap B_R(p)\su \bigcup_{x\in\3}B_{r_x}(x) \quad\text{with}\quad  r_x\leq r \quad\text{and}\quad  \sum_{x\in\3}r_x^k\leq C(m)R^k.
				\end{equation}
				Moreover, the balls in $\{B_{r_x/5}(x)\}_{x\in\cal{C}}$ are pairwise disjoint and $\cal{C}\subset S\cap B_R(p) $. 
			\end{lemma}
			
			\subsection{Proof of Lemma \ref{lemma: main covering lemma}}
			
			To prove Lemma \ref{lemma: main covering lemma}, we need two auxiliary lemmas. 
			
			\begin{lemma}[Covering Lemma I]\label{lemma: covering lemma I}
				{Let $J\in W^{1,2}_\La(\cJ_g(B_3(p)))$ be a minimizing harmonic almost complex structure on $B_3(p)\subset M$.} For any $\ep>0$, $0<\rho=\rho(m)\leq100^{-1}$ and $0<r<R\leq1$, there exist constants $\de=\de(m,g,\La,\rho,\ep)>0$ and $C_1(m)$ with the following property: 
				
				For any ${S}\su S^k_{\ep,\de r}$ with
				$$E:=\sup\limits_{x\in B_{2R}(p)\cap S}\Phi_J(x,R)\leq\Lambda, $$ 
				there exists a finite covering of $S\cap B_R(p)$ such that
				\begin{equation}\label{eq: covering lemma I}
					S\cap B_R(p)\su \bigcup_{x\in\3}B_{r_x}(x), \quad\text{with}\quad  r_x\geq r \quad\text{and}\quad  \sum_{x\in\3}r_x^k\leq C_1(m)R^k.
				\end{equation}
				Moreover, for each $x\in\3$, one of the following conditions is satisfied:
				\begin{itemize}
					\item[(i)] $r_x=r$;
					\item[(ii)] the set of points $F_x:=\{y\in S\cap B_{2r_x}(x):\Phi_J(y,\rho r_x/10)>E-\de\}$ is contained in $B_{\rho r_x/5}(L_x)\cap B_{2r_x}(x),$ where $L_x$ is some  $(k-1)$-dimensional affine subspace.
				\end{itemize}
			\end{lemma}


			\begin{lemma}[Covering Lemma II]\label{lemma: covering lemma II}
				{Let $J\in W^{1,2}_\La(\cJ_g(B_3(p)))$ be a minimizing harmonic almost complex structure on $B_3(p)\subset M$.} For any $\ep>0$ and $0<r<R\leq1$, there exist $\de=\de(m,g,\La,\ep)>0$ and $C_2(m)$ such that for any subset ${S}\su S^k_{\ep,\de r}$, there exists a finite covering of $S\cap B_R(p)$ by
				\begin{equation}\label{eq: covering lemma II-1}
					S\cap B_R(p)\su \bigcup_{x\in\3}B_{r_x}(x), \quad\text{with}\quad \ r_x\geq r\quad\text{and}\quad\sum_{x\in\3}r_x^k\leq C_2(m)R^k.
				\end{equation}
				Moreover, for each $x\in\3$, one of the following conditions is satisfied:
				\begin{itemize}
					\item[i)] $r_x=r;$
					\item[ii)] we have the following uniform energy drop property:
					\begin{equation}\label{eq: covering lemma II-2}
						\sup_{y\in B_{2r_x}(x)\cap S}\Phi_J(y,r_x)\leq E-\de,
					\end{equation} 
					where $$E:=\sup_{x\in B_{2R}(p)\cap{S}}\Phi_J(x,R).$$
				\end{itemize}
			\end{lemma}

			Based on the above two lemmas, we can complete the proof of Lemma \ref{lemma: main covering lemma}.
			\begin{proof}[Proof of Lemma \ref{lemma: main covering lemma}]
				Note that the energy $E$ defined as in Lemma \ref{lemma: covering lemma II} satisfies $E\leq C\Lambda$. So iterating Lemma \ref{lemma: covering lemma II} by at most $i=[\de^{-1}E]+1$ times, we could obtain a covering $\{B_{r_x}(x)\}_{x\in \3^i}$ of $S\cap B_R(p)$ such that $r_x\le r$ and 
				$$
				\sum_{x\in \cal{C}^i}r_x^k\leq C_2(m)^iR^k.
				$$ 
				Then for $x\in S\cap B_R(p)$, we have the larger covering $\{B_{2r_x}(x)\}_{x\in \tilde{\cal C}^i}$ such that
				$$
				S\cap B_R(p) \subset \bigcup_{x\in \tilde{\cal C}^i} B_{2 r_x}(x)\quad \text{ and } \quad  \sum_{x\in \cal{C}}(2 r_x)^k\leq 2^kC_2(m)^iR^k.
				$$
				Finally, by Vitali's covering lemma, we can  select a family of disjoint balls $\{B_{2r_x}(x)\}_{x\in \cal{C}}$ from  $\tilde{\cal C}^i$ such that
				$$ S\cap B_R(p) \subset \bigcup_{x\in \cal{C}} B_{10r_x}(x)\quad \text{ and } \quad  \sum_{x\in \cal{C}}(10r_x)^k\leq 10^kC_2(m)^iR^k.$$
				The proof is complete by taking $C(m)=10^kC_2(m)^i$.
			\end{proof}

			\subsection{Proof of Lemma \ref{lemma: covering lemma I}}
			
			{From now on, we denote $B_r(0)$ for $\exp_p^{-1}(B_r(p))$, and $S^k_{\e,r}$ for $\exp^{-1}\left(S^k_{\e,r}\right)$ for brevity.}

			\begin{proof}[Proof of Lemma \ref{lemma: covering lemma I}]
				We assume $R=1$. For the convenience, we also assume in the following proof that
				\begin{align}\label{eq: coverlemma-convenient notation}
					r=\rho^{{l}}, \quad \rho=2^{-a},\quad a,\ {l}\in\N.
				\end{align}
				
				We divide the proof into two main steps.
				\medskip 
				
				\noindent\textbf{Step 1. (Inductive covering)}
				\medskip 
				
				Let $\eta>0$ be a constant to be determined later. For a ball $B_r(z)$, if $\hat{F}:=\{y\in B_{2r}(z):\Phi_{J}(y,\rho/10)\ge E-\delta\}$ does not $\frac{\rho}{10}$-effectively span a $k$-dimensional affine subspace, then $\hat{F}$ is contained in $B_{\frac{\rho}{5}}(L)\cap B_{2r}(z)$ for some $(k-1)$-dimensional affine subspace $L$. We call $B_r(z)$ a \emph{good ball}, otherwise we call it a \emph{bad ball}.
				
				If $B_1(0)$ is a \emph{good ball}, then our claim clearly holds. 
				
				If $B_1(0)$ is a \emph{bad ball}, we define the initial set
				$$
				F^0=\{y\in B_2(0)\cap{S}:\Phi_J(y,\rho/10)\geq E-\de\}.
				$$
				By Proposition \ref{prop: stratum property-1}, if $\de\leq\de_2$ is sufficiently small, we  obtain
				\begin{equation}
					S^k_{\ep,\de r}\cap B_1(0)\su S^k_{\ep,\de}\cap B_1(0)\su B_{\rho/5}(V_0)
				\end{equation}
				for some $k$-dimensional affine subspace $V_0$, which is ${\rho}/{10}$-effectively spanned by some $\{y_j\}_{j=0}^k\su F^0$. Then we can choose a finite cover of $B_{\rho/5}(V_0)\cap B_1$ by balls $\{B_\rho(x)\}_{x\in\3^1}$ with $\3^1\su V_0\cap B_1(0)$ such that if $x\neq y$ and $x,y\in\3^1$, then $B_{\rho/5}(x)\cap B_{\rho/5}(y)=\emptyset.$
				
				Note that by Lemma \ref{lemma: Phi_J almost const}, if $\de\le\de_4$ sufficiently small, we have for all $x\in\3^1\su V_0\cap B_1(0)$,
				$$\Phi_J(x,\rho/10)\geq E-\eta,$$
				Under the smallness assumption for $\delta\le\de_5$, Lemma \ref{lemma: pinching almost symmetry} (or Remark \ref{rmk:on Lemma 313}) implies that for each $x\in \3^1$ we have $x\in S^k_{\ep/2,\rho}.$	
				Then we can divide $\3^1$ into two parts. Let $$
				\3^1_g:=\left\{x\in\3^1:\ F_x^1 \su B_{ \rho^2/5}(L_x^1) \
				\mbox{\rm for some}\ (k-1)\mbox{\rm -dimensional affine subspace}\ L_x^1\right\}
				$$and$$
				\3^1_b:=\{x\in\3^1:\ F_x^1\ (\rho^2/10)\mbox{-}\mbox{\rm effectively spans a}\ k\mbox{\rm -dimensional affine subspace}\ L_x^{1'}\},
				$$where$$
				F_x^1=\left\{y\in S\cap B_{2\rho}(x):\Phi_J(y, \rho^2/10)\geq E-\de\right\}.
				$$
				Then $$S\cap B_1(0)\subset\bigcup_{x\in \3^1}B_{\rho}(x)=\bigcup_{x\in \3^1_b}B_{\rho}(x)\cup\bigcup_{\3^1_g}B_{\rho}(x)\equiv B_{\rho}(\3^1_b)\bigcup B_{\rho}(\3^1_g).$$

				For any \emph{bad ball} $B_\rho(x)$. By Proposition \ref{prop: stratum property-1}, we have 
				$$S\cap B_\rho(x)\su B_{ \rho^2/5}(L_x^{1'})$$ for some $k$-dimensional affine subspace $L_x^{1'}$.
				Similarly, we can cover $B_{ \rho^2/5}(L_x^{1'})$ by balls $\{B_{\rho^2}(x)\}_{x\in\3^2_x}$ with $\3^2_x=\3^2_{x,b}\cup\3^2_{x,g}\su L_x^{1'}\cap B_\rho(x)$ so that if $y\neq z$ and $y,z\in\3^2_x$, then $B_{\rho^2/5}(y)\cap B_{\rho^2/5}(z)=\emptyset.$
				Write as above, $$S\cap B_\rho(x)\su\bigcup_{y\in \3^2_x}B_{\rho}(y)=\bigcup_{y\in \3^2_{x,b}}B_{\rho^2}(y)\cup\bigcup_{y\in \3^2_{x,g}}B_{\rho^2}(y)\equiv B_{\rho^2}(\3^2_{x,b})\bigcup B_{\rho^2}(\3^2_{x,g}).$$
				Here
				$$
				\3^2_{x,g}:=\{y\in\3^2_{x}:\ F_y^2 \su B_{ \rho^3/5}(L_y^2) \
				\ \ \mbox{\rm for some}\ (k-1)\mbox{\rm -dimensional affine subspace}\ L_y^2\}
				$$
				and
				$$
				\3^2_{x,b}:=\{y\in\3^2_x:\ F_y^2\ (\rho^3/10)\mbox{-}\mbox{\rm effectively spans a}\ k\mbox{\rm -dimensional affine subspace}\ L_y^{2'}\}.
				$$with
				$$
				F_y^2=\{z\in S\cap B_{2\rho}(y):\Phi_J(z, \rho^3/10)\geq E-\de\}.
				$$
				
				Let $\3_b^2=\bigcup\limits_{x\in\3^1_b}\3^2_{x,b}$, $\3^2_g=(\bigcup\limits_{x\in\3^1_b}\3^2_{x,g})\cup\3^1_g$ and $\3^2=\3^2_g\cup\3^2_b$. Then
				$$S\cap B_1(0)\su\bigcup_{x\in \3^2}B_{r^2_{x}}(x)=\bigcup_{x\in \3^2_g}B_{r^2_{x}}(x)\cup\bigcup_{x\in\3^2_b}B_{r^2_{x}}(x)\equiv B_{r^2_{x}}(\3^2_g)\bigcup B_{\rho^2}(\3^2_b),$$
				where $r^2_{x}=\rho$ if $x\in\3^1_g$ and $r^2_{x}=\rho^2$ if $x\in\3^2_g\ba\3^1_g$ or $\3^2_b$.
				Furthermore, with the same reasoning as the initial covering, there holds for all $x\in\3^2$ that $\Phi_J(x,r^2_{x})\geq E-\eta$ and for all $s\in[r^2_{x},1]$, $J$ is not $(k+1,\ep/2)$-symmetric on $B_{s}(x)$, that is, $x\in S^k_{\ep/2,r^2_{x}}.$
				
				After iterating, we then build a covering of the form
				$$
				S\cap B_1(0)\su \bigcup_{x\in \3^j}B_{r^j_{x}}(x)=\bigcup_{x\in \3^j_g}B_{r^j_{x}}(x)\cup\bigcup_{\3^j_b}B_{r^j_{x}}(x)\equiv B_{r^j_{x}}(\3^j_g)\bigcup B_{r^j_{x}}(\3^j_b),
				$$
				with $\3_b^j=\bigcup\limits_{x\in\3^{j-1}_b}\3^j_{x,b}$, $\3^j_g=\left(\bigcup\limits_{x\in\3^{j-1}_b}\3^j_{x,g}\right)\cup\3^{j-1}_g$ and $\3^j=\3^j_b\cup\3^j_g,j\geq2,$ where
				$$
				r^j_{x}=\begin{cases}
					\rho^i, & \text{ if } x\in\3^i_g\ba\3^{i-1}_g,i=2,\cdots,j,\\
					\rho^j, & \text{ if } x\in\3^j_b.
				\end{cases}
				$$
				Moreover, $\3^j_{g}$ is defined to be
				$$
				\{x\in\3^j:\ r^j_{x}\geq\rho^j\ \mbox{\rm and}\ F_x^j\su B_{ \rho r^j_{x}/5}(L_x^j) \
				\mbox{\rm for some}\ (k-1)\mbox{\rm -dimensional affine subspace}\ L_x^j\}
				$$
				and$$
				\3^j_{b}:=\{x\in\3^j: r^j_{x}=\rho^j \ \mbox{\rm and}
				\ F_x^j\ \frac{\rho r^j_{x}}{10}\mbox{-}\mbox{\rm effectively spans a}\ k\mbox{\rm -dimensional affine subspace}\ L_x^{j'}\}.
				$$with
				$$
				F_x^j:=\{y\in S\cap B_{2r^j_{x}}(x):\Phi_J(y, \rho r^j_{x}/10)\geq E-\de\}.
				$$
				Furthermore, we have
				\begin{itemize}
					\item[(P1)] For all $x\neq y\in\3^j\backslash \3_g^{j-1}$, $B_{r_x^j/5}(x)\cap B_{r_x^j/5}(y)=\emptyset$,
					\item[(P2)] For all $x\in\3^j,$  $\Phi_J(x,r^j_{x})\geq E-\eta.$ 
					\item[(P3)] For all $x\in\3^j$ and for all $s\in[r^j_{x},1]$, $J$ is not $(k+1,\ep)$-symmetric on $B_{s}(x)$, that is, $x\in S^k_{\ep/2,r^j_{x}}.$
				\end{itemize}
				
				Taking $j=l$ and $\rho^{l}=r$ in \eqref{eq: coverlemma-convenient notation} and $\3=\3^{l}$. Then the covering part of Lemma \ref{lemma: covering lemma I} is complete after selecting a subcollection of disjoint balls $\{B_{r_x}(x)\}$ such that $\{B_{5r_x}(x)\}$ satisfies all the covering requirements of the lemma. 
				\medskip 
				
				\noindent\textbf{Step 2. (Reifenberg estimates)}
				\medskip 
				
				In order to prove the volume estimate in \eqref{eq: covering lemma I}, we define a measure
				$$
				\mu:=\omega_k\sum_{x\in\3}r_x^k\de_x.
				$$
				and  measures
				$$
				\mu_t:=\omega_k\sum_{x\in\3_t}r_x^k\de_x,\qquad \text{with}\quad \3_t=\{x\in\3:r_x\leq t\},
				$$
				for all $t\in(0,1]$.

				Set $r_j=2^jr,j=0,1,\cdots,al-3$. Then $r_{al-3}=1/8$ by \eqref{eq: coverlemma-convenient notation}.
				Since $B_{r/5}(x)\cap B_{r/5}(y)=\emptyset$ for all $x\not=y\in\3_r,$ it follows easily that $\mu_r(B_r(x))\leq c(m)r^k.$
				Assume for all $x\in B_3(0)$ and $s\geq r,$ we have
				\begin{align}\label{eq: cover lemma I step2 induction-1}
					\mu_{r_j}(B_{r_j}(x)):=\left(\omega_k\sum_{x\in\3,r_x\leq r_j}r_x^k\de_x\right)\left(B_{r_j}(x)\right)\leq C_R(m)r_j^k,
				\end{align}
				where $C_R(m)$ is the constant in Theorem \ref{thm: discrete-Reifenberg}.  
				
				We first show that \eqref{eq: cover lemma I step2 induction-1} holds with constant $C_1(m)=c(m)C_R(m)$. For the convenience, we write
				$$
				\mu_{r_{j+1}}=\mu_{r_j}+\widetilde{\mu}_{r_{j+1}}:=\sum_{x\in\3_{{r_j}}}\omega_kr_x^k\de_x+\sum_{x\in\3,r_x\in({r_j},{r_{j+1}}]}\omega_kr_x^k\de_x.
				$$
				Take a covering of $B_{r_{j+1}}(x)$ by $M$ balls $\{B_{{r_j}}(y_i)\}$, $M\leq c(m)$, such that $\{B_{{r_j}/5}(y_i)\}$ are disjoint. Then we have
				$$
				\mu_{{r_j}}(B_{r_{j+1}}(x))\leq\sum_{j=1}^{M}\mu_{{r_j}}(B_{{r_j}}(y_i))\leq c(m)C_R(m)r_j^k.
				$$
				By definition of $\widetilde{\mu}_{r_{j+1}}$ and pairwise disjointness of $\{B_{r_x/5}(x)\}$, we have
				\[
				\widetilde{\mu}_{r_{j+1}}(B_{r_{j+1}}(x))\leq c(m)r_{j+1}^k.
				\]
				Thus for all $x\in B_1(0)$, there holds
				\begin{align}\label{eq: cover lemma I induction-2}
					\mu_{r_{j+1}}(B_{r_{j+1}}(x))\leq c(m)C_R(m)r_{j+1}^k.
				\end{align}
				This shows \eqref{eq: cover lemma I step2 induction-1} is true for ${j+1}$ and hence it holds for all $j\leq al-3$ by induction.
				
				For a fixed ball $B_{r_{j+1}}(x_0)$, set
				$$
				{\mu_{j+1}}:=\mu_{r_{j+1}}|_{B_{{r_{j+1}}}(x_0)},
				$$
				and we claim that for all $z\in \supp(\mu_{j+1})$
				\begin{align}\label{eq: cover lemma I Jones number claim}
					\beta_2=\beta^k_{2,{\mu_{j+1}}}(z,s)^2\leq C_1s^{-k}\int_{B_s(z)}\hat{W}_s(y)d{\mu_{j+1}}(y),
				\end{align}
				where $$
				\hat{W}_s(y):=\begin{cases}
					W_s(y), & \text{ if }s>r_z/5\\
					0, & \text{ if }0<s\le r_z/5.
				\end{cases}
				$$
				
				When $s\le r_z/5$, $$ \beta_2=\beta^k_{2,\mu_{j+1}}(z,s)^2=\inf_{V\su\R^m}\int_{B_s(z)}\frac{d^2(y,V)}{s^2}\frac{d\mu_{j+1}(y)}{s^k}=0$$
				since we can select $V$ passing $B_s(z)\cap \supp(\mu_{j+1})$, which is a single point.
				
				When $s>r_z/5$, since $\Phi_J(z,10s)\leq E$, by (P2) we have 
				$$\Phi_J(z,10s)-\Phi_J(z,5s)\leq\eta\quad \text{ for all } z\in \mbox{supp}(\mu) \text{ and all }s\in[ r_{z}/5,1/10].$$
				Given $\delta_8$ as in Theorem \ref{thm: L2 subspace approximation}, Proposition \ref{prop: quantitative symmetry} implies that there exists $\eta_0>0$ such that if $\eta\leq\eta_0$, then $J$ is $(0,\de_8)$-symmetric on $B_{10s}(x)$. By (P3), $J$ is not $(k+1,\ep/2)$-symmetric on  $B_{10s}(x)$. Then we have \eqref{eq: cover lemma I Jones number claim} holds for $0<s\leq 1/10$ by applying Theorem \ref{thm: L2 subspace approximation}.
				
				Note that, by the induction assumption and \eqref{eq: cover lemma I induction-2}, for all $j\leq al-3$ and $s\in(0,r_{j+1}],$ and $z\in B_1(0),$ we have
				\begin{equation}\label{eq:rough growth for mu s}
					\mu_s(B_s(z))\leq c(m)C_R(m)s^k.
				\end{equation}
				We claim that for any $r<s\leq r_{j+1}$, we have
				\begin{equation}\label{eq:improved growth for mu j}
					\mu_{r_{j+1}}(B_s(z))\leq c(m)C_R(m)5^ks^k.
				\end{equation}
				Indeed, if $y\in B_s(z)\cap \supp(\mu)$, then $\frac{r_y}{5}\leq |y-z|\leq s$ and so $y\in \3_{5s}$, which implies $B_s(z)\cap \supp(\mu)\subset \3_{5s}$. Since $r\leq 5s\leq 5r_{j+1}$, we have
				\[
				\mu_{j+1}(B_s(z))\leq \mu_{5s}(B_s(z))\leq \mu_{5s}(B_{5s}(z))\leq c(m)C_R(m)5^ks^k.
				\]
				
				Thus, 
				\begin{equation}\label{eq: cover lemma I final cal}
					\begin{aligned}
						&\quad\int_{B_r(y)}\left(\int_0^r\beta^k_{2,\mu_{j+1}}(z,s)^2\frac{ds}{s}\right)d\mu_{j+1}(z)=  \int_0^r\frac{ds}{s}\int_{B_r(y)}\beta^k_{2,\mu_{j+1}}(z,s)^2d\mu_{j+1}(z)\\
						&\leq  C_1\int_0^r\frac{ds}{s}\int_{B_r(y)}s^{-k}\left(\int_{B_s(z)}\hat{W}_s(y)d{\mu_{j+1}}(y)\right)d\mu_{j+1}(z)\\
						&\leq  C_1\int_0^r\frac{ds}{s^{k+1}}\int_{B_r(y)}\int_{B_{2r}(y)}\chi_{B_s(z)}(x)W_s(x)d\mu_{j+1}(x)d\mu_{j+1}(z)\\
						&=  C_1\int_0^r\frac{ds}{s^{k+1}}\int_{B_{2r}(y)}\mu_{j+1}(B_s(z))d\mu_{j+1}(x)\\
						&\leq  C_1c(m)C_R(m)\int_0^r\frac{ds}{s}\int_{B_{2r}(y)}W_s(x)d\mu_{j+1}(x)\\
						&= C_1c(m)C_R(m)\int_{B_{2r}(y)}d\mu_{j+1}(x)\int_0^rW_s(x)\frac{ds}{s}\\
						&\le C_1c(m)C_R(m)C\eta r^k.
					\end{aligned}
				\end{equation}
				Here in the third inequality wee used \eqref{eq:improved growth for mu j} and in the last inequality, we used$$
				\int_0^r{W}_s(x)\frac{ds}{s}=\int_{r_x/5}^r{W}_s(x)\frac{ds}{s}\leq\int_{r_x/5}^{1/10}{W}_s(z)\frac{ds}{s}\leq C[\Phi_J(x,1)-\Phi_J(x,r_x/5)]\leq C\eta
				$$for $x\in\supp(\mu)$ and $r\leq r_{j+1}\leq1/8$.

				Then the estimate \eqref{eq: cover lemma I Jones number claim} follows from Theorem \ref{thm: discrete-Reifenberg} applied to $\mu_{j+1}$ when choosing $\eta$ small enough.
			\end{proof}

			\subsection{Proof of Lemma \ref{lemma: covering lemma II}}
			
			\begin{proof}[Proof of Lemma \ref{lemma: covering lemma II}]
				For simplicity we assume that $R=1$. We will use an induction argument to refine the first covering lemma so as to deduce the desired covering. The proof is divided into several steps. We will use superscripts $f,b$ to indicate \emph{final} and \emph{bad} balls respectively. We will also call the  properties i) and ii) in Lemma \ref{lemma: covering lemma II} as stopping conditions in our proof. 
				\medskip 
				
				\noindent\textbf{Step 1. (Initial covering)} 
				\medskip 
				
				By the previous lemma, we have a covering of $S\cap B_{1}(0)$ given by
				$$
				S\cap B_{1}(0)\subset\bigcup_{x\in{\cal C}_{r}^{0}}B_{r}(x)\cup\bigcup_{x\in{\cal C}_{+}^{0}}B_{r_{x}}(x),
				$$
				where
				$$
				{\cal C}_{r}^{0}=\{x\in{\cal C}:r_{x}=r\}\quad\text{and}\quad{\cal C}_{+}^{0}=\{x\in{\cal C}:r_{x}>r\},
				$$
				and
				\begin{equation}\label{eq: cover lemma II-1}
					\sum_{x\in{\cal C}_{r}^{0}\cup{\cal C}_{+}^{0}}\om_{k}r_{x}^{k}\le C_1(m)=:C_{V}(m).
				\end{equation}
				Moreover, for each $x\in{\cal C}_{+}^{0}$, the set
				\[
				F_{x}=\left\{ y\in S\cap B_{2r_{x}}(x):\Phi_J(y, \rho r_{x}/10)>E-\de\right\} 
				\]
				is contained in a small neighborhood of a ($k-1$)-dimensional affine subspace. We only need to refine the part ${\cal C}_{+}^{0}$. 
				\medskip 
				
				\noindent\textbf{Step 2. (Recovering of bad balls: the first step)} 
				\medskip 
				
				Let $x\in{\cal C}_{+}^{0}$. If $\rho r_{x}=r$, then we just cover $S\cap B_{r_{x}}(x)$
				by a family of balls $\{B_{r_{y}}(y)\}_{y\in{\cal C}_{x}^{(1,r)}}$
				with $\{B_{\rho r_{x}/2}(y)\}_{y\in{\cal C}_{x}^{(1,r)}}$ being pairwise disjoint, where $y\in B_{r_{x}}(x)$, $r_{y}=\rho r_{x}=r$.  A simple volume comparison argument shows that $\sharp\left({\cal C}_{x}^{(1,r)}\right)\le C(m)\rho^{-n}$, which gives
				$$
				\sum_{y\in{\cal C}_{x}^{(1,r)}}r_{y}^{k}=\sharp\left({\cal C}_{x}^{(1,r)}\right)(\rho r_{x})^{k}\le C(m)\rho^{k-m}r_{x}^{k}=:C_{r}(m,\rho)r_{x}^{k}.
				$$
				Thus we collect all such points $x$ and set
				$$
				{\cal C}^{(1,r)}:={\cal C}_{r}^{0}\cup\bigcup_{x\in{\cal C}_{+}^{0},\rho r_{x}=r}{\cal C}_{x}^{(1,r)}.
				$$
				By \eqref{eq: cover lemma II-1}, we have
				$$
				\sum_{y\in{\cal C}^{(1,r)}}r_{y}^{k}=\left(\sum_{x\in{\cal C}_{r}^{0}}+\sum_{x\in{\cal C}_{+}^{0}}\sum_{y\in{\cal C}_{x}^{(1,r)}}\right)r^{k}\le\sum_{x\in{\cal C}_{r}^{0}}r^{k}+C_{r}(m,\rho)\sum_{x\in{\cal C}_{+}^{0}}r_{x}^{k}\le 2C_{V}(m)C_{r}(m,\rho).
				$$
				
				Next, suppose $\rho r_{x}>r$. Two cases occur.
				
				\smallskip 
				\noindent\textbf{Case 1: $F_{x}=\emptyset$.}
				\smallskip 
				
				In this case we simply cover $S\cap B_{r_{x}}(x)$
				by balls $\{B_{r_{y}}(y)\}_{y\in{\cal C}_{x}^{(1,f)}}$ centered in
				$S\cap B_{r_{x}}(x)$ with $r_{y}=\rho r_{x}$, such that $\{B_{\rho r_{x}/2}(y)\}_{y\in{\cal C}_{x}^{(1,f)}}$
				are disjoint. In this case, the energy drop property clearly holds and the number $\sharp\left\{ {\cal C}_{x}^{(1,f)}\right\} $
				is bounded from above by a constant $C(m)\rho^{-m}$, so that
				$$
				\sum_{y\in{\cal C}_{x}^{(1,f)}}r_{y}^{k}=\sharp\left({\cal C}_{x}^{(1,f)}\right)(\rho r_{x})^{k}\le C(m)\rho^{k-m}r_{x}^{k}=:C_{f}(m,\rho)r_{x}^{k}.
				$$
				Note that actually we have $C_{f}(m,\rho)=C_{r}(m,\rho)$. 
				\smallskip 
				
				\noindent\textbf{Case 2: $F_{x}\neq\emptyset$.}
				\smallskip 
				
				We call $B_{r_{x}}(x)$ a bad ball. We want to recover $B_{r_{x}}(x)$ using the fact that $F_{x}\subset B_{\rho r_{x}/5}(L_{x}^{k-1})\cap B_{2r_{x}}(x)$.
				This will be done as follows. First note that the part away from $B_{\rho r_{x}}(F_{x})$
				is good: we have
				\[
				S\cap B_{r_{x}}(x)\Big\backslash B_{\rho r_{x}}(F_{x})\subset\bigcup_{y\in{\cal C}_{x}^{(1,f)}}B_{r_{y}}(y)\quad\text{with }r_{y}=\rho r_{x},
				\]
				and $\{B_{\rho r_{x}/2}(y)\}_{{\cal C}_{x}^{(1,f)}}$ being pairwise
				disjoint. Since $d(y,F_x)\ge\rho r_x>r$, the energy drop property holds:
				$\Phi_J(y,r_{y}/10)\le E-\de$ for all $y\in{\cal C}_{x}^{(1,f)}$.
				Thus these balls will be part of the ``final'' balls in the next
				step. Moreover, we also have the trivial estimate:
				$$
				\sum_{y\in{\cal C}_{x}^{(1,f)}}r_{y}^{k}=\sharp\left\{ {\cal C}_{x}^{(1,f)}\right\} (\rho r_{x})^{k}\le C(m)\rho^{k-n}r_{x}^{k}=C_{f}(m,\rho)r_{x}^{k}.
				$$
				Then we collect all the ${\cal C}_{x}^{(1,f)}$ from the above to
				obtain final balls of the first generation:
				$$
				{\cal C}^{(1,f)}=\bigcup_{x\in{\cal C}_{+}^{0},\rho r_{x}>r}{\cal C}_{x}^{(1,f)}.
				$$
				By \eqref{eq: cover lemma II-1}, we have
				$$
				\sum_{y\in{\cal C}^{(1,f)}}r_{y}^{k}=\sum_{x\in{\cal C}_{+}^{0}}\sum_{y\in{\cal C}_{x}^{(1,f)}}(\rho r_{x})^{k}\le C_{f}(m,\rho)\sum_{x\in{\cal C}_{+}^{0}}r_{x}^{k}\le C_{V}(m)C_{f}(m,\rho).
				$$
				
				For the remaining part of $S\cap B_{r_{x}}(x)$, we simply
				cover it by
				$$
				S\cap B_{r_{x}}(x)\cap B_{\rho r_{x}}(F_{x})\subset\bigcup_{y\in{\cal C}_{x}^{1,b}}B_{r_{y}}(y)\quad\text{with }r_{y}=\rho r_{x},
				$$
				and $\{B_{\rho r_{x}/2}(y)\}_{{\cal C}_{x}^{(1,b)}}$ being pairwise
				disjoint, where "$b$" means bad balls. The problem is that the
				energy drop condition can not be verified on each ball $B_{r_{y}}(y)$.
				However, since $F_{x}\subset B_{\rho r_{x}/5}(L_{x})\cap B_{2r_{x}}(x)$
				for some ($k-1$)-dimensional affine subspace $L_{x}$, a volume comparison
				argument gives 
				$$
				\sharp\{{\cal C}_{x}^{(1,b)}\}\le C(m)\frac{(\rho r_x/5)^{m-k+1}(2r_x)^{k-1}}{(\rho r_x)^m}\le C(m)\rho^{1-k},
				$$
				This means that there are relatively very few bad balls. Hence
				$$
				\sum_{y\in{\cal C}_{x}^{(1,b)}}(\rho r_{x})^{k}=(\rho r_{x})^{k}\sharp\left\{ \text{\ensuremath{{\cal C}_{x}^{(1,b)}}}\right\} \le C(m)\rho r_{x}^{k}=:C_{b}(m)\rho r_{x}^{k}.
				$$
				Thus we collect all the bad balls and define
				$$
				{\cal C}^{(1,b)}=\bigcup_{x\in{\cal C}_{+}^{0},\rho r_{x}>r}{\cal C}_{x}^{(1,b)}.
				$$ 
				
				Now we have the following properties for balls in ${\cal C}^{(1,b)}$
				and ${\cal C}^{(1,r)}\cup{\cal C}^{(1,f)}$. For bad balls,  using \eqref{eq: cover lemma II-1}, we obtain
				\[
				\sum_{y\in{\cal C}^{(1,b)}}r_{y}^{k}=\sum_{x\in{\cal C}_{+}^{0}}\sum_{y\in{\cal C}_{x}^{(1,b)}}(\rho r_{x})^{k}\le C_{b}(m)\rho\sum_{x\in{\cal C}_{+}^{0}}r_{x}^{k}\le C_{V}(m)C_{b}(m)\rho.
				\]
				Here and after, we choose
				$$
				0<\rho<\min\left\{ 100^{-1},\frac{1}{2C_{V}(m)C_{b}(m)}\right\} 
				$$
				such that
				$$
				\sum_{y\in{\cal C}^{(1,b)}}(\rho r_{x})^{k}<1/2.
				$$
				Set
				$$
				C_{2}=C_{2}(m)=2C_{V}(m)\left(C_{r}(m,\rho)+C_{f}(m,\rho)\right),
				$$
				then we have
				$$
				\sum_{y\in{\cal C}^{(1,r)}\cup{\cal C}^{(1,f)}}r_{y}^{k} \le(C_{f}(m,\rho)+C_{r}(m,\rho))\sum_{x\in{\cal C}}r_{x}^{k}\le\frac{1}{2}C_{2}(m).
				$$
				
				Note that for all $y\in{\cal C}^{(1,b)}$, $r<r_{y}\le\rho$. 
				\medskip

				\noindent\textbf{Step 3. (Induction step)} 
				\medskip 
				
				We claim the existence of following covering: for each $i\ge1$,
				$$
				S\cap B_{1}(0)\subset\bigcup_{x\in{\cal C}^{(i,r)}}B_{r}(x)\cup\bigcup_{x\in{\cal C}^{(i,f)}}B_{r_{x}}(x)\cup\bigcup_{x\in{\cal C}^{(i,b)}}B_{r_{x}}(x)
				$$
				with the following properties hold:
				
				(1) For each $x\in{\cal C}^{(i,r)}$, $r_{x}=r$;
				
				(2) For each $x\in{\cal C}^{(i,f)}$, the energy drop condition holds:
				for all $z\in \cal{S}\cap B_{2r_{x}}(x)$ we have $\Phi_J(z,r_{x}/10)\le E-\de$;
				these balls will be called final balls;
				
				(3) For each $x\in{\cal C}^{(i,b)}$, we have $r<r_{x}\le\rho^{i}$. On these "bad" balls, none of the above two stopping conditions is verified. 
				
				(4) There holds
				\begin{equation}\label{eq: cover lemma II final estimate}
					\sum_{y\in{\cal C}^{(i,r)}\cup{\cal C}^{(i,f)}}r_{y}^{k}\le\left(\sum_{j=1}^{i}2^{-j}\right)C_{2}(m)\quad\text{and}\quad \sum_{y\in{\cal C}^{(i,b)}}r_{y}^{k}\le2^{-i}.
				\end{equation}

				We have proved the above induction for $i=1$ in Step 2. Suppose now it holds for some $i\ge1$. We aim to prove that it holds for $i+1$. Clearly, we only need to recover these bad balls. 
				
				Suppose now $x\in{\cal C}^{(1,b)}$, i.e., $B_{r_{x}}(x)$ is a bad ball and $r<r_{x}\le\rho^{i}$.
				\smallskip 
				
				\noindent\textbf{Case 1: $\rho r_{x}=r$.}
				\smallskip 
				
				We cover every $B_{r_x}(x)$ simply by $\{B_{\rho r_x}(y)\}_{y\in\3^{(i+1,r)}_x}$ with $\{B_{\rho r_x/2}(y)\}$ pairwise disjoint as in step 2. Then we have the estimate
				$$
				\sum_{y\in{\cal C}_{x}^{(i+1,r)}.}(\rho r_{x})^{k}\le C(m)\rho^{k-n}r_{x}^{k}=C_{r}(m,\rho)r_{x}^{k}.
				$$

				\smallskip 
				\noindent\textbf{Case 2: $\rho r_{x}>r$.} 
				\smallskip 
				
				In this case we apply Lemma \ref{lemma: covering lemma I} to get a covering $\{B_{r_{y}}(y)\}_{y\in{\cal C}^{x}}$ such that
				$$
				S\cap B_{r_{x}}(x)\subset\bigcup_{y\in{\cal C}_{r}^{x}}B_{r}(y)\cup\bigcup_{y\in{\cal C}_{+}^{x}}B_{r_{y}}(y)\quad\text{with }r_{y}\ge r
				$$
				and
				$$
				\sum_{y\in{\cal C}_{r}^{x}\cup{\cal C}_{+}^{x}}r_{y}^{k}\le C_{V}(m)r_{x}^{k}.
				$$
				Moreover, for each $y\in{\cal C}_{+}^{x}$, there is a ($k-1$)-dimensional
				affine subspace $L_{y}$ such that
				$$
				F_{y}\equiv\left\{ z\in S\cap B_{2r_{y}}(y):\Phi_J(z,\rho r_{y}/10)\ge E-\de\right\} \subset B_{\rho r_{y}/5}(L_{y})\cap B_{2r_{y}}(y).
				$$
				Below we assume that $y\in{\cal C}_{r}^{x}$. 
				\smallskip 
				
				\noindent\textbf{Case 2.1: $\rho r_{y}=r$.}
				\smallskip 
				
				As that of Case 1, we get a simple covering of at most $C(m)\rho^{-m}$ balls $\{B_{r}(z)\}_{z\in{\cal C}_{y}^{(i+1,r)}}$ of $B_{r_{y}}(y)$. Now we define
				$$
				{\cal C}^{(i+1,r)}={\cal C}^{(i,r)}\cup\bigcup_{x\in{\cal C}^{(i,b)},\rho r_{x}=r}{\cal C}_{x}^{(i+1,r)}\cup\bigcup_{x\in{\cal C}^{(i,b)},\rho r_{x}>r}\left({\cal C}_{r}^{x}\cup\bigcup_{y\in{\cal C}_{+}^{x},\rho r_{y}=r}{\cal C}_{y}^{(i+1,r)}\right).
				$$
				This shows that how much more is ${\cal C}^{(i+1,r)}$ than that of
				${\cal C}^{(i,r)}$.
				\smallskip 
				
				\noindent\textbf{Case 2.2: $\rho r_{y}>r$.}
				\smallskip 
				
				If $F_y=\emptyset$, we cover $S\cap B_{r_{y}}(y)$ by$$
				S\cap B_{r_{y}}(y)\subset\bigcup_{z\in{\cal C}_{y}^{(i+1,f)}}B_{r_{z}}(z)\quad\text{with }r_{z}=\rho r_{y}.
				$$

				If $F_y\not=\emptyset$, we cover $B_{r_{y}}(y)$ by
				$$
				\begin{aligned}
					& S\cap B_{r_{y}}(y)\Big\backslash B_{\rho r_{y}}(F_{y})\subset\bigcup_{z\in{\cal C}_{y}^{(i+1,f)}}B_{r_{z}}(z),\\
					& S\cap B_{r_{y}}(y)\cap B_{\rho r_{y}}(F_{y})\subset\bigcup_{z\in{\cal C}_{y}^{(i+1,b)}}B_{r_{z}}(z),
				\end{aligned}
				\quad\text{with }r_{z}=\rho r_{y}.
				$$
				Here each ball $B_{r_{z}}(z)$ with $z\in{\cal C}_{y}^{(i+1,f)}$ satisfies the energy drop condition. Also we have the following estimates:
				$$
				\sum_{z\in{\cal C}_{y}^{(i+1,f)}}r_{z}^{k}\le C_{f}(m,\rho)r_{y}^{k}\quad\text{and}\quad \sum_{z\in{\cal C}_{y}^{(i+1,b)}}r_{z}^{k}\le C_{b}(m)\rho r_{y}^{k}.
				$$
				We define
				$$
				{\cal C}^{(i+1,f)}={\cal C}^{(i,f)}\cup\bigcup_{x\in{\cal C}^{(i,b)}}\bigcup_{y\in{\cal C}_{+}^{x}}{\cal C}_{y}^{(i+1,f)},\quad{\cal C}^{(i+1,b)}=\bigcup_{x\in{\cal C}^{(i,b)}}\bigcup_{y\in{\cal C}_{+}^{x}}{\cal C}_{y}^{(i+1,b)}.
				$$

				Now we prove \eqref{eq: cover lemma II final estimate}. For $\cal C^{(i+1,r)}$, we have
				\begin{equation}\label{eq: cover lemma II Cr estimate}
					\begin{aligned}
						\sum_{z\in{\cal C}^{(i+1,r)}}r_{z}^{k} & =\sum_{z\in{\cal C}^{(i,r)}}r_{z}^{k}+\sum_{x\in{\cal C}^{(i,b)},\rho r_{x}=r}\sum_{z\in{\cal C}_{x}^{(i+1,r)}}r_{z}^{k}\\
						& \quad+\sum_{x\in{\cal C}^{(i,b)},\rho r_{x}>r}\left(\sum_{z\in{\cal C}_{r}^{x}}r_{z}^{k}+\sum_{y\in{\cal C}_{+}^{x}}\sum_{z\in{\cal C}_{y}^{(i+1,r)}}r_{z}^{k}\right)\\
						& \le\sum_{z\in{\cal C}^{(i,r)}}r_{z}^{k}+C_{r}(m,\rho)\left(\sum_{x\in{\cal C}^{(i,b)},\rho r_{x}=r}r_{x}^{k}+\sum_{x\in{\cal C}^{(i,b)},\rho r_{x}>r}\sum_{y\in{\cal C}^{x}}r_{y}^{k}\right)\\
						& \le \sum_{z\in{\cal C}^{(i,r)}}r_{z}^{k}+C_{r}(m,\rho)C_{1}(m)\sum_{x\in{\cal C}^{(i,b)}}r_{x}^{k}\\
						& \le\sum_{z\in{\cal C}^{(i,r)}}r_{z}^{k}+2^{-i}C_{r}(m,\rho)C_{1}(m).
					\end{aligned}
				\end{equation}
				For $\cal C^{(i+1,f)}$, we have
				\begin{equation}\label{eq: cover lemma II Cf estimate}
					\begin{aligned}
						\sum_{z\in{\cal C}^{(i+1,f)}}r_{z}^{k} & =\sum_{z\in{\cal C}^{(i,f)}}r_{z}^{k}+\sum_{x\in{\cal C}^{(i,b)}}\sum_{y\in{\cal C}_{+}^{x}}\sum_{z\in{\cal C}_{y}^{(i+1,f)}}r_{z}^{k}\\
						& \le\sum_{z\in{\cal C}^{(i,f)}}r_{z}^{k}+C_{f}(m,\rho)\sum_{x\in{\cal C}^{(i,b)}}\sum_{y\in{\cal C}_{+}^{x}}r_{y}^{k}\\
						& \le\sum_{z\in{\cal C}^{(i,f)}}r_{z}^{k}+C_{f}(m,\rho)C_{V}(m)\sum_{x\in{\cal C}^{(i,b)}}r_{x}^{k}\\
						& \le\sum_{z\in{\cal C}^{(i,f)}}r_{z}^{k}+2^{-i}C_{f}(m,\rho)C_{V}(m).
					\end{aligned}
				\end{equation}
				Recall that $C_{2}(m)=2\left(C_{r}(m,\rho)+C_{f}(m,\rho)\right)C_{V}(m)$, combining \eqref{eq: cover lemma II Cr estimate} and \eqref{eq: cover lemma II Cf estimate} we have 
				$$
				\left(\sum_{z\in{\cal C}^{(i+1,r)}}+\sum_{z\in{\cal C}^{(i+1,f)}}\right)r_{z}^{k}\le\left(\sum_{z\in{\cal C}^{(i,r)}}+\sum_{z\in{\cal C}^{(i,f)}}\right)r_{z}^{k}+2^{-i-1}C_{2}(m)\le C_{2}(m)\sum_{j=1}^{i+1}2^{-j}.
				$$
				Also by the choice of $\rho$, we obtain
				$$
				\begin{aligned}\sum_{z\in{\cal C}^{(i+1,b)}}r_{z}^{k}= & \sum_{x\in{\cal C}^{(i,b)}}\sum_{y\in{\cal C}_{+}^{x}}\sum_{z\in{\cal C}_{y}^{(i+1,b)}}r_{z}^{k}\le\sum_{x\in{\cal C}^{(i,b)}}\sum_{y\in{\cal C}_{+}^{x}}C_{b}(m)\rho r_{y}^{k}\\
					& \le\sum_{x\in{\cal C}^{(i,b)}}C_{b}(m)\rho C_{1}(m)r_{x}^{k}\le2^{-1}\sum_{x\in{\cal C}^{(i,b)}}r_{x}^{k}\le2^{-i-1}.
				\end{aligned}
				$$
				This proves \eqref{eq: cover lemma II final estimate} for $i+1$.
				
				Note that for every $z\in{\cal C}^{(i+1,b)}$, we have $r_{z}=\rho r_{y}\le\rho r_{x}\le \rho^{i+1}$ for some $x\in{\cal C}^{(i,b)}$.  The proof of Step 3 is complete.

				Since $r=\rho^{\bar{j}}$, the above procedure will stop at $i=\bar{j}-1$. 
				\medskip 
				
				\noindent\textbf{Step 4. (Final refinement) }
				\medskip 
				
				Now we have obtained a covering of $S\cap B_R(0)$ satisfying 
				$$
				\begin{aligned}
					S\cap B_R(0)\su \bigcup_{x\in\3}B_{r_x}(x)=\bigcup_{x\in\3_{r}}B_{r_x}(x)\cup\bigcup_{x\in\3_{+}}B_{r_x}(x),\quad  \sum_{x\in\3}r_x^k\leq C_2(m),
				\end{aligned}
				$$
				where $\3_{r}$ consists of centers of all balls $B_{r_x}(x)$ such that $r_x=r$ and  $\3_{+}$ consists of centers of all balls $B_{r_x}(x)$ such that $r_x> r$ and
				$$
				\begin{aligned}
					\sup_{y\in S\cap B_{2r_x}(x)}\Phi_J({y,r_x/10})\leq E-\de.
				\end{aligned}
				$$
				For each $x\in\3_{+}$, we can simply cover  $S\cap B_{r_x}(x)$ by a family of balls $\{B_{r_y}(y)\}_{y\in \3_x}$  such that $y\in S\cap B_{r_x}(x) $,  $\{B_{r_y/5}(y)\}_{y\in \3_x}$ are pairwise disjoint and $\sharp\{\3_x\}\leq C(m)\rho^{-m}$ with $r_y=\rho r_x$. Then we have $B_{2r_y}(y)\subset B_{2r_x}(x)$ and $r_y\le r_x/10$, which implies by monotonicity that
				$$
				\sup_{z\in S\cap B_{2r_y}(y)}\Phi_J(z,r_y)\leq \sup_{z\in S\cap B_{2r_x}(x)}\Phi_J(z,r_x/10)\leq E-\delta. 
				$$
				Moreover, 
				$$
				\sum_{x\in \3_+}\sum_{y\in \3_x}r_y^k\leq C(m)\rho^{k-m}\sum_{x\in \3_+}r_x^k\leq C(m)\rho^{k-m}C_2(m).
				$$
				This completes the proof of Lemma \ref{lemma: covering lemma II} by taking $C_2'(m)=2C(m)\rho^{k-m}C_2(m)$.
			\end{proof}

			\section{Proofs of main theorems}\label{sec:proofs}
			
			\subsection{Proof of Theorem \ref{thm: stratification of MHACS}}
			Theorem \ref{thm: stratification of MHACS} follows from the main covering lemma, Lemma \ref{lemma: main covering lemma}, and the rectifiable Reifenberg theorems in Section \ref{sec:reifenberg}.
			\begin{proof}[Proof of Theorem \ref{thm: stratification of MHACS}]
				
				Let $\de=\de(m,\Lambda,\ep)$ be fixed. By Lemma \ref{lemma: main covering lemma}, for all $0<r<\delta$,
				\begin{equation}\label{eq: main thm proof-1}
					\mbox{Vol}(T_r(S^k_{\ep,r}(J))\cap B_1(p))\le \mbox{Vol}\Big(T_r\Big(S^k_{\ep,\de }(J)\cap B_1(p)\Big)\Big) \leq C'_\ep r^{m-k}
				\end{equation}
				with $C'_{\ep}=C'_\ep(m,\Lambda,\ep)$. For $r\in[\de,1)$, we have 
				\begin{equation}\label{eq: main thm proof-2}
					\mbox{Vol}\Big(T_r\Big(S^k_{\ep,r}(J)\cap B_1(p)\Big)\Big) \le \mbox{Vol}(T_1( B_1(p)) \leq C_m \le C_m \left(\frac{r}{\de} \right)^{m-k}.
				\end{equation}
				\eqref{eq: main thm proof-1} and \eqref{eq: main thm proof-2} complete the proof of  \eqref{eq: Vol estomate-1} for all $0<r<1$.  The volume estimate \eqref{eq: Vol estomate-2} follows from \eqref{eq: Vol estomate-1} by noting that $S^k_{\ep}(J)\subset S^k_{\ep, r}(J)$ for any $r>0$.

				In order to prove the rectifiability of $S^k(J)$, it is sufficient to prove the rectifiability of $S^k_\ep(J)$ for each $\ep>0$, as $S^k(J)=\bigcup\limits_{i\ge 1}S^k_{1/i}(J)$.
				
				By the volume estimate \eqref{eq: Vol estomate-2}, we have $\HH^k(S^k_\ep(J)\cap B_1(p))\leq C_\ep.$ Similarly, for $B_r(x)$ with $x\in B_1(p)$ and $r\leq 1$, we have
				\begin{equation}\label{eq: main thm Alhfors}
					\HH^k(S^k_\ep(J)\cap B_r(x))\leq C_\ep r^k.
				\end{equation}
				
				Let $S\su S_\ep^k(J)\cap B_1(p)$ be an arbitrary measurable subset with $\HH^k(S)>0$. Set 
				$$g(x,r)=\Phi_J(x,r)-\Phi_J(x,0),\quad x\in B_1(p),\ r\leq 1.$$
				This function is nondecreasing in $r$ and uniformly bounded for all $x\in B_1(p)$ and $r\in[0,1]$, and pointwise converging to $0$ as $r\to0$. Then the dominated convergence theorem implies that for all $\de>0$, there exists $\overline{r}>0$ such that
				\begin{equation}\label{eq: select E}
					\medint_S g(x,10\overline{r})d\HH^k(x)\leq\de^2.
				\end{equation}
				So we can find a measurable subset $E\su S$ with $\HH^k(E)\leq\de\HH^k(S)$ and $g(x,10\overline{r})\leq\de$ for all $x\in F:=S\ba E$. Now Cover $F$ by a finite number of balls $B_{\overline{r}}(x_i)$ centered on $F$. Rescaling if necesssary, we may assume that $B_{\overline{r}}(x_i)=B_1(0)$. Then $g(x,10)\leq\de$ for $x\in F$. Similar to \eqref{eq: select E}, we can choose $\de$ sufficiently small such that $J$ is $(0,\de_9)$ symmetric in $B_{10}$. Theorem \ref{thm: L2 subspace approximation} implies 
				\begin{equation*}
					\beta_{2,\HH^k|_F}(z,s)^2\leq C_1s^{-k}\int_{B_s(z)}W_s(t)d\HH^k|_F(t) \quad \text{for all }z\in F,\ s\leq1.
				\end{equation*}
				Integrating the previous estimate with respect to $z$ and using \eqref{eq: main thm Alhfors} yield that for all  $x\in B_1(p)$ and $s\leq r\leq1,$ there holds
				\begin{align*}
					\int_{B_r(x)}\beta_{2,\HH^k|_F}(z,s)^2d\HH^k|_F(z)&\leq C_1s^{-k}\int_{B_r(x)}\int_{B_s(z)}W_s(t)d\HH^k|_F(t)d\HH^k|_F(z)\non\\
					&\leq C_1C_\ep\int_{B_{r+s}(x)}W_s(z)d\HH^k|_F(z).
				\end{align*}
				Integrating again with respect to $s$, similar to \eqref{eq: cover lemma I final cal}, we obtain that for all $x\in B_1(p)$ and $r\leq1,$
				\begin{align*}
					\int_{B_r(x)}\int_0^r\beta_{2,\HH^k|_F}(z,s)^2\frac{ds}{s}d\HH^k|_F(z)&\leq C_1C_\ep\int_{B_{2r}(x)}[\Phi_J(z,10r)-\Phi_J(z,0)]d\HH^k|_F(z)\non\\
					&\leq c(m)C_1C_\ep^2\de r^k.
				\end{align*}
				By Choosing $\de\le\frac{\de_7^2}{c(m)C_1C_\ep^2}$,
				we deduce from Theorem \ref{thm: rectifiable-Reifenberg} that $F\cap B_1(p)$ is $k$-rectifiable.
				
				Repeating the above argument with $E$ in place of $F$, we could find another measurable set $E_1\subset E$ with $\HH^k(E_1)\leq \delta \HH^k(E)$, and that $F_1:=E\ba E_1$ is $k$-rectifiable. Continuing this process, we eventually conclude that $S$ is $k$-rectifiable. 
				
				The proofs of final assertions are similar to that of \cite{Naber-V-2017} and are thus omitted.
			\end{proof}

			\subsection{Proof of Theorem \ref{thm: regularity estimates on minimizing HACS}}
			
			The proof of Theorem \ref{thm: regularity estimates on minimizing HACS} is built on the important compactness results. First of all, we prove the following theorem. 
			
			\begin{theorem}[Symmetry implies regularity] \label{thm: HACS symmetry implies regularity}
				There exists $\de(m,g,\La)>0$ such that if $J\in W^{1,2}_\La(\cJ_g(B_1(p)))$ is an $(m-2,\de)$-symmetric minimizing harmonic almost complex structure on $B_1(p)\subset M$, then
				$$
				r_{J}(p)\ge1/2.
				$$
			\end{theorem}
			
			The proof is divided into the following two lemmas. 
			
			\begin{lemma} [{$(m,\ep)$-regularity theory}]\label{lemma:new epsilon regularity} 
				There exists $\ep=\ep(m,g,\La)>0$ such that if $J\in W^{1,2}_\La(\cJ_g(B_1(p)))$ is an $(m,\ep)$-symmetric minimizing harmonic almost complex structure on $B_1(p)\subset M$, then
				$$r_{J}(p)\ge1/2.$$
			\end{lemma}
			
			Note that the $(m,\ep)$-regularity theory is an improvement of the usual partial regularity Theorem \ref{thm: ep-regularity}.
			
			\begin{proof}
				Suppose by contradiction that there is a sequence of minimizing harmonic almost complex structure $J_{i}\in W^{1,2}_{\La}(\cJ_g)$ satisfying
				$$
				\medint_{B_{1}}|J_{i}-C_{i}|^{2}<1/i
				$$for a sequence of constants $\{C_i\}$, but $J_{i}\not\in C^{\infty}(B_{1/2})$. Up to subsequence if necessary, we may assume that $J_{i}\wto J_{\infty}$ weakly in $W^{1,2}$ and $C_i\to C_\infty$. By Theorem \ref{thm: compact minimizing}, we have $J_{i}\to J_{\infty}$ strongly in $W^{1,2}$. Sending $i\to\infty$, we find that $J_{\infty}$ is a constant map. Then we have
				$$
				\medint_{B_{1}}|\na J_{i}|^{2}\to\medint_{B_{1}}|\na J_{\infty}|^{2}=0.
				$$
				Thus, the $\ep$-regularity theorem \ref{thm: ep-regularity} implies $J_{i}\in C^{\infty}(B_{1/2})$ for $i\gg1$. This is a contradiction.
			\end{proof}
			
			The second lemma gives us improvement of symmetry.
			\begin{lemma}[Symmetry improvement]\label{lemma:symmetry improvement}
				For any $\ep>0$, there exists $\de=\de(m,g,\La,\ep)>0$ such that if $J\in W^{1,2}_\La(\cJ_g(B_1(p)))$ is an $(m-2,\de)$-symmetric minimizing harmonic almost complex structure on $B_1(p)\subset M$, then $J$ is $(m,\ep)$-symmetric on $B_{1}(p)$.
			\end{lemma}
			\begin{proof}
				Suppose that there is a constant $\ep>0$ such that, for each $i\ge1$ there is a minimizing harmonic almost complex structure $J_i\in W^{1,2}_\La(\cJ_g(B_1))$, which is $(m-2,1/i)$-symmetric on $B_{1}$ but not $(m,\ep)$-symmetric on $B_{1}$. Up to a subsequence if necessary, we have $J_{i}\wto J_{\infty}$ weakly in $W^{1,2}$. By the compactness Theorem \ref{thm: compact minimizing}, we know that $J_{\infty}$ is a minimizing $(m-2)$-symmetric harmonic almost complex structure, but not $(m,\ep/2)$-symmetric on $B_{1}$.
				
				Thus $J_{\infty}$ can be viewed as a $0$-homogeneous minimizing  harmonic almost complex structure from $B_1^2$ to $TM\otimes T^*M$. Hence $J_{\wq}$ is smooth by Remark \ref{rmk: 2-dim HACS is smooth}. Thus, the $0$-homogeneity of $J_{\wq}$ implies that $J_\infty$ is a constant. This leads to a contradiction.
				
				
				
			\end{proof}

			\begin{proof}[Proof of Theorem \ref{thm: HACS symmetry implies regularity}]
				It follows straightforward from Lemmas \ref{lemma:new epsilon regularity} and \ref{lemma:symmetry improvement}.
			\end{proof}

			Now we sketch the proof of Theorem \ref{thm: regularity estimates on minimizing HACS}. 
			\begin{proof}[Proof of Theorem \ref{thm: regularity estimates on minimizing HACS}]
				The First estimate of \eqref{eq: nabla J in L3,infty} is directly from Definition \ref{def: regularity scale}, which implies
				\begin{align*}
					\{x\in B_1(p):r|\nabla J|+r^2|\nabla^2 J|>1\}\su\{x\in B_1(p):r_J(x)<r\}.
				\end{align*}

				By a scaling argument, Theorem \ref{thm: HACS symmetry implies regularity} implies
				$$\{x\in B_1(p) : r_J(x)<r\}\su S_{\ep,r}^{m-3}(J).$$
				Thus by Theorem \ref{thm: stratification of MHACS} there exists $C>0$ such that for each $0<r<1$ we have
				$$
				{\rm Vol}(T_r(\{x\in B_1(p):r_J(x)<r\}))\leq{\rm Vol}(T_r(S_{\ep,r}^{m-3}(J)))\leq Cr^3,
				$$
				which gives the second estimate of \eqref{eq: nabla J in L3,infty}. The proof is thus complete. 
			\end{proof}
			
			\appendix
			
			\section{Proof of Theorem \ref{thm: L2 subspace approximation} }
			
			Without loss of generality, we assume that $\mu$ is a probability measure supported on $B_1(0)$. Let $x_{cm}$ be the mass center of $\mu$ in $B_1(0)$, i.e.,  
			$
			x_{cm}\equiv\int xd\mu(x).
			$
			Here and later the integral domain is  $B_1(0)$ and is omitted for brevity. The second moment $Q=Q(\mu)$ of $\mu$ is the symmetric bilinear form defined by
			$$
			Q(v,w):=\int[(x-x_{cm})\cdot v][(x-x_{cm})\cdot w]d\mu(x),\qquad \text{for all }\ v,w\in \R^m.
			$$
			Let $\lambda_1(\mu)\ge \cdots\ge \lam_m(\mu)$ be nonincreasing eigenvalues of $Q(\mu)$ and $v_1(\mu),\cdots,v_m(\mu)$ be the associated eigenvectors. Then we have
			\begin{align}\label{eq: Q eigenvalue}
				Q(v_k)=\lam_kv_k=\int[(x-x_{cm})\cdot v_k][(x-x_{cm}) ]d\mu(x).
			\end{align}By the variational method, we can describe $\lam_1$ by $$\lam_1=\lam_1(\mu):=\max_{|v|^2=1}\int|(x-x_{cm})\cdot v|^2d\mu(x),$$
			where $v_1=v_1(\mu)$ is any unit vector achieving such maximum. By induction we have
			$$
			\lam_{k+1}=\lam_{k+1}(\mu):=\max\left\{\int|(x-x_{cm})\cdot v|^2d\mu(x):{|v|^2=1}, v\perp\mbox{span}\{v_i\}, \text{for all }i\leq k\right\},
			$$
			and $v_{k+1}=v_{k+1}(\mu)$ is any unit vector achieving such maximum. Note that by definition of $v_k$, it is not difficult to show that  $V_k=x_{cm}+\mbox{span}\{v_1,\cdots,v_k\}$ is the $k$-dimensional affine subspace achieving the minimum in the definition of $\beta_2$ (see \cite[Remark 49]{Naber-V-2018}). Moreover,
			$$
			\beta_{2,\mu}^k(0,1)^2=\int d^2(x,V_k)d\mu(x)=\int\sum_{i=k+1}^{m}((x-x_{cm})\cdot v_i)^2d\mu=\lam_{k+1}(\mu)+\cdots+\lam_m(\mu).
			$$
			To prove Theorem \ref{thm: L2 subspace approximation}, we need the following property about the relationship of $\lam_k,v_k$ and $W_1$. Here and after, let $A_{r_1,r_2}(a):=B_{r_2}(a)\backslash B_{r_1}(a)$ and $A_{r_1,r_2}:=B_{r_2}(p)\backslash B_{r_1}(p)$.
			
			\begin{proposition}[{\cite[Proposition 6.3]{Naber-V-2018}}]\label{prop: relationship of lam,v,W}
				Let $\mu$ be a probability measure in $B_1(0)$ and a map $J\in H^1(B_{10}(0), N)$ for some compact manifold $N$. Let $\lam_k$, $v_k$ be defined as above. Then there exists $C(m)>0$ such that
				$$
				\lam_k\int_{A_{3,4}}|v_k\cdot\na J(z)|^2dz\leq C(m)\int_{B_1} W_1(x)d\mu(x), \qquad \text{for all }\,k\ge 1,
				$$
			\end{proposition}
			
			
			For convenience of the readers, we sketch the proof here.
			\begin{proof}[Proof of Proposition \ref{prop: relationship of lam,v,W}] 
				Without loss of generality, we assume  $x_{cm}=0$. For any $z\in A_{3,4}$ and $k=1,\cdots,m$, multiplying $\na J(z)$ on both sides of \eqref{eq: Q eigenvalue} yields
				\begin{align}\label{eq: Prop 4.5-1} 		\lam_k(v_k\cdot\na J(z))=\int(x\cdot v_k)(\na J(z)\cdot x)d\mu(x).	
				\end{align}
				By definition of mass center, for all $z$ we have
				$$
				\int x\cdot zd\mu(x)=x_{cm}\cdot z=0.
				$$
				Hence by  \eqref{eq: Prop 4.5-1}, 	
				\begin{align*}	
					\lam_k(v_k\cdot\na J(z))=\int(x\cdot v_k)(\na J(z)\cdot (x-z))d\mu(x).
				\end{align*}
				By H\"{o}lder inequality we have
				\begin{align*}		
					\lam_k^2|v_k\cdot\na J(z)|^2\leq\lam_k\int |(\na J(z)\cdot (x-z))|^2d\mu(x).	\end{align*}
				Noticing we can assume $\lam_k>0,$ otherwise there is nothing to prove. Partial integration leads to
				\begin{align*}		\lam_k\int_{A_{3,4}(0)}|\na J(z)\cdot v_k|^2dz&\leq \int\int_{A_{3,4}(0)}|\na J(z)\cdot(x-z)|^2dzd\mu(x)\\		
					&\leq \int\int_{A_{3,4}(0)}\frac{|\na J(z)\cdot(x-z)|^2}{|x-z|^m}|x-z|^mdzd\mu(x)\\
					&\leq C(m) \int\int_{A_{1,10}(x)}\frac{|\na J(z)\cdot(x-z)|^2}{|x-z|^m} dzd\mu(x)\\
					&\leq C(m)\int W_1(x)d\mu(x).	
				\end{align*}
				This completes the proof.
			\end{proof}
			
			To further estimate the left hand side of the estimate in the above proposition, we need the following reversed result of Lemma \ref{lemma: almost translation invariant}.
			\begin{lemma}\label{lemma: non higher order sym}
				For any $\ep>0$, there exists  $\de_9=\de_9(m,g,\La,\ep)>0$ satisfying the following property. For any minimizing harmonic almost complex structure {$J\in W^{1,2}_\La(\cJ_g(B_2(p)))$, if it is $(0,\de_9)$-symmetric on $B_1(p)$} but not $(k+1,\ep)$-symmetric, then 	
				\begin{equation}\label{eq: non higher order sym}
					\int_{A_{3,4}}|P\cdot D J|^2>\delta_9
				\end{equation}
				for every $(k+1)$-dimensional subspace $P$.
			\end{lemma}
			\begin{proof}
				Assume by contradiction that there is a sequence of minimizing harmonic almost complex structures {$J_i\in W^{1,2}_\La(\cJ_g(B_2(p)))$ such that $J_i$ is $(0,i^{-1})$-symmetric on $B_1(p)$} but not $(k+1,\ep_0)$-symmetric for some $\ep_0>0$ fixed. Moreover, after an orthogonal transformation, there is a  $(k+1)$-dimensional subspace $P$ such that
				$$
				\int_{A_{3,4}}|P\cdot D J_i|^2\leq i^{-1},\ i=1,2,\cdots.
				$$
				By Theorem \ref{thm: compact minimizing}, up to a subsequence if necessary, we may assume $J_i\to J$ strongly in $W^{1,2}$ for some $0$-homogeneous minimizing harmonic almost complex structure $J$ and 
				\begin{align*}
					\int_{A_{3,4}}|P\cdot D J|^2=0.
				\end{align*}
				By the unique continuation property Proposition \ref{prop: unique continuation}, $\int_{B_{4}}|P\cdot D J|^2=0.$ Thus $J$ is $(k+1)$-symmetric on $B_1(p)$. Consequently $J_i$ is $(k+1,\ep_0)$-symmetric for $i\gg 1$, since $J_i$ converges strongly to $J$ in $L^2$. This is a contradiction. 
			\end{proof}
			
			Now we can prove Theorem \ref{thm: L2 subspace approximation}.
			\begin{proof}[Proof of Theorem \ref{thm: L2 subspace approximation}]  
				We denote $B_1(0)$ for $\exp_p^{-1}(B_1(p))$. Without loss of generality, we assume $\mu(B_1(0))=1$. Since $\lambda_k$ is nonincreasing,
				\begin{equation}\label{eq: beta<(m-k)lam-k+1}
					\beta_{2,\mu}^k(0,1)^2=\lam_{k+1}+\cdots+\lam_m\leq(m-k)\lam_{k+1}.
				\end{equation}
				Thus it suffices to estimate $\la_{k+1}$. 
				Let $V^{k+1}=\mbox{span}(v_1,\cdots,v_{k+1}),$ then
				\begin{align*}
					&\lam_{k+1}\int_{A_{3,4}}|V^{k+1}\cdot DJ(z)|^2dz=\lam_{k+1}\sum_{j=1}^{k+1}\int_{A_{3,4}}|DJ(z)\cdot v_j|^2dz\\
					\le& \sum_{j=1}^{k+1}\lam_j\int_{A_{3,4}}|DJ(z)\cdot v_j|^2dz\le C\int W_1d\mu(x),
				\end{align*}
				where the second inequality is due to Proposition \ref{prop: relationship of lam,v,W}.
				Set $\delta_8=\delta_9$, where $\delta_9$ is the number defined in Lemma \ref{lemma: non higher order sym}. By Lemma \ref{lemma: non higher order sym}, there holds 
				\begin{align*}
					\int_{A_{3,4}}|V^{k+1}\cdot DJ(z)|^2dz\geq \de_8,
				\end{align*}
				since by assumption $J$  is $(0,\de_8)$-symmetric but not $(k+1,\ep)$-symmetric on $B_{10}(p)$. Then we have 
				\begin{align*}
					\de_8\lam_{k+1}\leq\lam_{k+1}\int_{A_{3,4}}|V^{k+1}\cdot DJ(z)|^2dz \leq C\int W_1d\mu(x).
				\end{align*}
				Since $\de_8=\de_9$ depends only on $m,g,\La,\ep$, by \eqref{eq: beta<(m-k)lam-k+1}, we conclude 
				\begin{align*}
					\beta_{2,\mu}^k(0,1)^2\leq C(m,\La,\ep)\int W_1d\mu(x).
				\end{align*}
				This completes the proof.
			\end{proof}

			\bigskip 
			\textbf{Acknowledgements.} 
			C.-Y. Guo and M.-L. Liu are supported by the Young Scientist Program of the Ministry of Science and Technology of China (No.~2021YFA1002200), the National Natural Science Foundation of China (No.~12101362), the Taishan Scholar Project and the Natural Science Foundation of Shandong Province (No.~ZR2022YQ01). C.-L. Xiang is supported by the NFSC grant (No.~12271296) and NSF of Hubei Province (No.~2024AFA061).
			
			We would like to thank Prof.~Weiyong He for his very valuable communications and insightful comments. In particular, his question on rectifiability and optimal regularity of harmonic almost complex structures leads to the current article. We also want to thank Prof.~Xiao Zhong for his helpful comments that improves our presentation.


\begin{thebibliography}{10}
				
				
				\bibitem{Aronszajn-1957}
				\textsc{N. Aronszajn}, \emph{A unique continuation theorem for solutions of elliptic partial differential equations or Inequalities}, J. Math. Pures Appl. (9), \textbf{36} (1957), 235–249.
				
				\bibitem{Calabi-Gluck-1993}
				\textsc{E. Calabi and H. Gluck}, \emph{What are the best almost-complex structures on the 6-sphere?} Proc. Sympos. Pure Math., 54, Part 2, Amer. Math. Soc., Providence, RI, 1993. 
				
	
				
				
				\bibitem{Cheeger-Naber-2013}
				\textsc{J. Cheeger and A. Naber}, \emph{Lower bounds on Ricci curvature and quantitative behavior of sigular sets}, Invent. Math. \textbf{191} (2013), 321-339.
				
				\bibitem{Cheeger-Naber-2013-CPAM}
				\textsc{J. Cheeger and A. Naber}, \emph{Qantitative stratification and the regularity of harmonic maps and minimal currents}, Comm. Pure. Appl. Math. \textbf{66} (2013), no.6, 965-990.
				
				\bibitem{Cheeger-N-V-2015}
				\textsc{J. Cheeger, A. Naber and D.Valtorta}, \emph{Critical sets of elliptic equations}, Comm. Pure Appl. Math. \textbf{68} (2015), no. 2, 173-209.
				
				\bibitem{Cheeger-J-N-2021}
				\textsc{J. Cheeger, W. Jiang and A. Naber},  \emph{Rectifiability of singular sets of noncollapsed limit spaces with Ricci curvature bounded below}, Ann. of Math. (2) \textbf{193} (2021), no. 2, 407-538.
				
				\bibitem{CLMS-1993} \textsc{R. Coifman, P.-L. Lions, Y. Meyer and S. Semmes,} \emph{Compensated compactness and Hardy spaces.} J. Math. Pures Appl. (9) \textbf{72} (1993),  247-286.
				
				\bibitem{Davidov-2017}
				\textsc{J. Davidov}, \emph{Harmonic almost Hermitian structures}, Special metrics and group actions in geometry, 129-159, Springer INdAM Ser., 23, Springer, Cham, 2017.
				
				\bibitem{Edelen-E-2019-TAMS}
				\textsc{N. Edelen and M. Engelstein}, \emph{Quantitative stratification for some free-boundary problems}, Trans. Amer. Math. Soc. \textbf{371} (2019), no. 3, 2043-2072.
				
				\bibitem{Evans-1991}\textsc{C.L. Evans,} \emph{ Partial regularity for stationary harmonic maps into spheres.} Arch. Rat. Mech. Anal. \textbf{116} (1991), 101-163.
				
				\bibitem{Fu-Wang-Zhang-2024}
				\textsc{H. Fu, W. Wang and Z. Zhang}, \emph{Quantitative stratification and sharp regularity estimates for supercritical semilinear elliptic equations}, arXiv-preprint, available at \url{https://arxiv.org/abs/2408.06726}, 2024. 
				
				\bibitem{Giaquinta-Book}\textsc{M. Giaquinta}, \emph{Multiple integrals in the calculus of variations and nonlinear elliptic systems.} Annals of Mathematics Studies, 105. Princeton University Press, Princeton, NJ, 1983.
				
				
				
				\bibitem{GJXZ-2024} \textsc{C.-Y. Guo, G.-C Jiang, C.-L. Xiang and G.-F. Zheng}, \emph{Optimal higher regularity for biharmonic maps via quantitative stratification}, Peking Math. J.(2025), available at \url{https://doi.org/10.1007/s42543-025-00107-0}.
				
				\bibitem{Guo-Liu-Xiang-2025}\textsc{C.-Y. Guo, M.-L. Liu and C.-L. Xiang}, \emph{$L^p$-regularity of a geometrically nonlinear system in supercritical dimensions}, Sci. China Math. \textbf{68} (2025), 2313-2332 

				
				\bibitem{He-2019} \textsc{W. Y. He}, \emph{Energy minimizing harmonic almost complex structures}, Ann. Mat. Pura Appl. (4), published online, available at \url{https://doi.org/10.1007/s10231-025-01632-6}, 2025.
				
				\bibitem{He-2023}
				\textsc{W. Y. He}, \emph{Biharmonic almost complex structures}, Forum Math. Sigma 11 (2023), Paper No. e25.
				
				\bibitem{He-Jiang-2021}
				\textsc{W. Y. He and R. Jiang}, \emph{Polyharmonic almost complex structures}, J. Geom. Anal. 31 (2021), no. 12, 11648-11684. 
				
				\bibitem{Helein-2002} \textsc{F. H\'elein,}
				\emph{Harmonic maps, conservation laws and moving frames.} Cambridge Tracts in Mathematics, \textbf{150}. Cambridge University Press, Cambridge, 2002.
				
				
				\bibitem{Lin-1999-Annals} \textsc{F. H. Lin}, \emph{Gradient estimates and blow-up analysis for stationary harmonic maps}, Ann. of Math. (2) \textbf{149} (1999), no. 3, 785-829.	
				
				\bibitem{Lin-Wang-2008-book}
				\textsc{F.H. Lin and C.Y. Wang}, \emph{The analysis of harmonic maps and their heat flows}, World Scientific Publishing Co. Pte. Ltd., Hackensack, NJ, 2008.
				
				\bibitem{Luckhaus-1988-Indiana} \textsc{S. Luckhaus}, \emph{Partial H\"older continuity for minima of certain energies among maps into a Riemannian manifold}. Indiana Univ. Math. J. 37 (1988), no. 2, 349-367.
				
				\bibitem{Naber-V-2017}
				\textsc{A. Naber and D. Valtorta}, \emph{Reifenberg-rectifiable and the regularity of stationary and minimizing harmonic maps}, Ann. of Math. (2) \textbf{185} (2017), no.1, 131-227.
				
				\bibitem{Naber-V-2018}
				\textsc{A. Naber and D. Valtorta}, \emph{Stratification for sigular set of approximate harmonic maps}, Math. Z. \textbf{290} (2018), no. 3-4, 1415-1455.
				
				\bibitem{Naber-V-2020-JEMS}
				\textsc{A. Naber and D. Valtorta}, \emph{The singular structure and regularity of stationary varifolds}, J. Eur. Math. Soc. (JEMS) \textbf{22} (2020), no. 10, 3305-3382.
				
				
				\bibitem{Riviere-2007-Invent}
				\textsc{T. Riviere}, \emph{Conservation laws for conformally invariant variational problems}, Invent. Math. \textbf{168} (2007), 1-22.
				
				
				\bibitem{Schikorra-2010} \textsc{A. Schikorra, }  \emph{A remark on gauge transformations and the moving frame method}. Ann. Inst. H. Poincar\'e C Anal. Non Lin\'eaire \textbf{27} (2010), no. 2, 503-515.
				
				\bibitem{SU-1982-JDG} \textsc{R. Schoen, K. Uhlenbeck}, \emph{A regularity theory for harmonic mappings.} J. Diff. Geom. 17 (1982), 307-335.
				
				\bibitem{Simon-1996-book} \textsc{L. Simon}, \emph{Theorems on regularity and singularity of energy minimizing maps}. Lectures in Mathematics ETH Zurich. Birkh\"auser Verlag, Basel, 1996.
				
				
				\bibitem{Wood-1993-Crelle}
				C.M. Wood, \emph{Instability of the nearly-K\"ahler six-sphere}, J. Reine Angew. Math. \textbf{439} (1993), 205-212.
				
				\bibitem{Wood-1995-Comp}
				C.M. Wood, \emph{Harmonic almost-complex structures}, Compositio Math. \textbf{99} (1995), no. 2, 183-212. 
				
				
							
				
			\end{thebibliography}
		\end{document}